\documentclass[letterpaper, 10pt]{article}

\usepackage{graphicx}

\usepackage{blindtext} 
\usepackage{amssymb,amsfonts,amsmath}

\usepackage[sc]{mathpazo} 
\usepackage[T1]{fontenc} 
\usepackage{microtype} 

\usepackage[english]{babel} 

\usepackage[letterpaper, hmargin=0.5in, vmargin=1in]{geometry}
\usepackage[hang, small,labelfont=bf,up,textfont=it,up]{caption} 
\usepackage{booktabs} 

\usepackage{lettrine} 
\usepackage{titlesec} 

\usepackage{enumitem} 
\setlist[itemize]{noitemsep} 

\usepackage{cite}

\usepackage{abstract} 

\renewcommand\thesection{\Roman{section}} 
\renewcommand\thesubsection{\roman{subsection}} 

\def\defi{:=}

\titleformat{\section}[block]{\large\scshape\centering}{\thesection.}{1em}{} 
\titleformat{\subsection}[block]{\large}{\thesubsection.}{1em}{} 

\usepackage{fancyhdr} 
\usepackage{titling} 

\usepackage{hyperref} 

\allowdisplaybreaks

\pretitle{\begin{center}\Huge\bfseries} 
\posttitle{\end{center}} 
\title{Joint Estimation of Discrete Choice Model and Arrival Rate with Unobserved Stock-out Events} 
\author{%
\textsc{Hongzhang Shao}\thanks{Georgia Institute of Technology (email: steveshao@gatech.edu)} \quad 
\textsc{Anton Kleywegt}\thanks{Georgia Institute of Technology (email: anton@isye.gatech.edu)} \\ [1ex] 
}
\date{\today} 

\renewcommand{\maketitlehookd}{%
\begin{abstract}

\noindent This paper studies the joint estimation problem of a discrete choice model and the arrival rate of potential customers when unobserved stock-out events occur. In this paper, we generalize \cite{anupindi1998estimation} and [Conlon and Mortimer, 2013] in the sense that (1) we work with generic choice models, (2) we allow arbitrary numbers of products and stock-out events, and (3) we consider the existence of the null alternative, and estimates the overall arrival rate of potential customers. In addition, we point out that the modeling in [Conlon and Mortimer, 2013] is problematic, and present the correct formulation. \\

\end{abstract}
}

\begin{document}

\maketitle

\section{Introduction}

Consider a vending machine that is replenished occasionally by a worker. At each visit to the vending machine, the worker first records the remaining amount of inventory of each product from the previous time period. After that, he determines the new assortment of products to be placed in the vending machine for the next time period, and replenish the inventory of each product to the desired stock level. (There are many more products than can be accommodated in a vending machine.) The updated initial assortment and stock levels for the next time period are also recorded. With these records, the sales of these products between any adjacent two replenishments can be determined. \\

Suppose that we want to estimate a demand model for this vending machine, in order to facilitate price and assortment optimization. To do that, we need to estimate the overall arrival rate of customers, as well as a discrete choice model that describes the choice behavior of the customers. Assume that each consumer chooses her favorite option $a$ from a choice set $\mathcal{A}$ provided, and makes a purchase if it is better than the no-purchase option (null alternative) $o$. That is, the choice probability $P_{a:\mathcal{A} \cup \{o\}}$ depends on the product being chosen, and the assortment being chosen from. (For simplicity in notations, external attributes are not considered in this paper, but can be accommodated easily if necessary.) \\ 

In this scenario, there are two challenges in estimating the demand model. First, when demand exceeds the available inventory, products go out-of stock, and the assortment becomes smaller. If we ignore the stock-out events, then the excess demand will not be captured, and the estimation can thus be serious biased. To make an unbiased estimation, it is preferred that we observe each individual choice and choice set. However, since only sales data is recorded, different customers may have different unobserved choice sets. If some products become out of stock between successive replenishment visits, we do not to know the order that products become out of stock, and thus do not know what choice sets are encountered. Even when there is only one stock-out event, we still do not know how many customers encounter each of the two choice sets. \\

Second, customers can choose to purchase nothing (the null alternative). However, we can only observe sales of products from data, and do not know the number of customers who choose the null alternative. If we assume that everyone makes a purchase, then the optimal prices should be astronomically high, which is certainly not realistic. The action of choosing the null alternative being unobserved makes it very hard to estimate the choice model and the overall arrival rate of potential customers. \\

In this paper, we study the estimation of generic demand models under different scenarios, by defining the underlying random process, and deriving the corresponding likelihood functions. To make things clear, we discuss three different types of data: complete data (penal data), transaction data and sales data. In complete data, we observe the entire sequence of choices and choice sets. Such data can usually be collected by online retailers based on browsing histories. In transaction data, we observe the entire sequence of transactions (true purchase), and the corresponding choice sets. Such data is commonly used by retailers like supermarkets, where any arrival that does not lead to a purchase is not observed. In sales data, the cumulative sales of products are observed periodically. Such data is more common with vending machines (like the one we just described) and traditional grocery stores. Here, the "sale" of the null alternative is not observed, and the order of the transactions is also unknown. Thus, although the initial assortments and inventory levels are known, we may still not be able to know the choice set associate with each choice. \\

We also want to point out that, whether or not the transaction times are observed makes a large difference to the likelihood of transaction data. This fact has been largely ignored, as most researchers just assume the time points to be known. Intuitively, as the assortment shrinks over time, the arrival rate of transactions also gets smaller. Thus, the cumulative arrival rate of transactions depends on the time when the assortment changes (stock-out times). On the contrary, observing the arrival times does not add new information to the estimation of the choice models when complete data is given. Meanwhile, in sales data, the transaction times must be unknown. (Otherwise, we are able to know the order of the transactions, and can then infer the choice sets. Such data cannot be sales data.) \\


\subsection{Background \& Related Literature}

Revenue management is the science to analyze consumer behavior with mathematical models, and to maximize revenue at the level of micro-market. The primary aim of revenue management is to choose the best product availability, price, pack, etc., and to sell it to the right customer at the right time. Businesses face important decisions regarding what to sell, when to sell, to whom to sell, and for how much. Revenue management answers these questions with data and models in order to increase profit. The concept of revenue management was mainly originated from the airline industry, who obtained the freedom of pricing its own products after the Civil Aviation Board in 1978. To survive the fierce competition, it soon became vital for airline companies to carefully manage their prices of the seats. Seeing the great value, other industries like ride-sharing, retailing, and many others, quickly adopted the revenue management technique as well. \\

Demand estimation and forecasting is the gate and the key to most of revenue manage practices, including assortment planning, inventory planning, pricing, etc.. Traditional revenue management assumes there is an independent stream of customers for each product, and use different demand models to describe what a customer may do when her preferred product is not available. However, after dozens of years, this assumption is becoming less and less realistic, even for the airline industry who first adopted it. In the 1980s, people usually book airline tickets through travel agents, who have to repeatedly sending requests to the booking system until available tickets are found. Now a customer can go online, and directly pick her favorite option from a large choice set. The substitution behavior of customers has thus been greatly changed. The limitations of the historical demand models and better data availability led to the adoption of discrete choice models in revenue management. \\ 

Before being adopted in revenue management, discrete choice models have already been widely used in explanatory (economics, marketing, transportation, etc.) studies for decades. The foundation of choice modeling comes from many early economical and psychological literatures. As an example, [Thurston, 1927] starts from the axiom that consumers make choices which maximize their perceived utility, and introduces the concept of random utility. It suggests that, since it is almost impossible to measure the utility of each individual customer, it is more practical to described the utility of the population with a random distribution. On the other hand, \cite{luce1959individual} introduces an axiom of choices, that for any set $\mathcal{A}' \subseteq \mathcal{A}$ such that $a \in \mathcal{A}'$, we have $P_{a:\mathcal{A}} = P_{a:\mathcal{A}'} \cdot P_{\mathcal{A}':\mathcal{A}}$. Choice models that follow this axiom have many properties. Among them, the irrelevant alternatives (IIA) property has been discussed in many later works. \\ 

[McFadden, 1973] combines the two ideas by showing that a random utility model with a Gumbel distribution is a conditional logit model that satisfies the choice axiom from \cite{luce1959individual}. Their model came to be known as the multinomial logit (MNL) model afterwards. By far, the MNL model is the most popular choice model in literature as well as in industrial usage, since both its estimation and optimization problems are tractable. However, many researchers believe that the irrelevant alternatives (IIA) property is too restrictive, and have been trying to get around it by introducing more generalized choice models. We refer to \cite{ben1985discrete} and \cite{train2009discrete} for detailed introductions of these choice models, including the nested logit (NL) model and the mixed logit (MMNL) model. In addition to those parametric models, there are also non-parametric choice models introduced in recent years, For example, \cite{blanchet2016markov} develops a Markov chain choice model that generalizes many earlier choice models. \\ 

There are many works that discuss the estimation of the attractiveness of the null alternative, and the total number of potential customers. For example, \cite{talluri2004revenue} studies the joint estimation problem of MNL model with the null alternative and the arrival rate. \cite{vulcano2012estimating} provides an efficient algorithm to solve this problem. In contrast, there are relatively less literature that talks about estimating choice models when stock-out events occur. \cite{anupindi1998estimation} is an early work along this line, which works with a generic demand model with two products. Another work, [Conlon and Mortimer, 2013], estimates the MNL model with stock-outs, which focus mainly on the case when there is only one stock-out event. Both works ignore the existence of the null-alternative. \\

We need to point out that the formulation in [Conlon and Mortimer, 2013] is not correct. [Conlon and Mortimer, 2013] suggests that the total sales (of other product) before the stock-out event follows a conditional negative binomial distribution. Given this total sales, the individual sales of each product before the stock-out event is distributed independently as a binomial random variable, whose parameter depends on that total sales. The problem here is that, since that total sales before the stock-out event is given, the distribution of the individual sales cannot be independent. As an example, suppose that there are only two products $1$ and $2$. Product $1$ is out of stock at some point, and $2$ is always available. In this case, that total sales is nothing but the individual sale of $2$ before $1$ is out of stock. In other words, the individual sale is fixed. We will provide a more concrete counter-example in a later section. \\

\subsection{Contribution \& Arrangement} 

This paper studies the joint estimation problem of discrete choice models and the total arrival rate of potential customers. The major contribution of this paper, is that we provide a general modeling framework to estimate choice models with stock-out events and the null alternative. Comparing to the earlier works, we work with generic choice models, and allow an arbitrary number of products and stock-out events. \\

We start this work by defining the underlining random process, and deriving the likelihood functions of complete data. This is presented in Section II. In Section III, we formulate likelihood functions of transaction data, and show that if the transaction times are unobserved, the likelihood function becomes much more complicated. In Section IV, we formulate the likelihood functions of sales data. We also discuss the form of the likelihood function under attraction demand models, and show that the size of the sums can be reduced with such choice models. In Section V, we discuss the likelihood function of sales data when the null alternative is assumed to be not available. Based on that discussion, we present a counter-example to show that the formulation in [Conlon and Mortimer, 2013] is not correct. In Section VI, we discuss the way we can deal with infinite sums, and introduce an approach to evaluate the likelihood functions by sampling. Numerical results are presented in Section IV. \\

\section{The Random Process}

Consider a store that offers products to customers. Assume that the store can only be in one of the two states: either customers arrive according to a Poisson process with a constant rate $\lambda$ (operating hours), or with rate $0$ (nights, weekends, holidays, etc.). Ignoring the idle time, we can model the customer arrival process as a stationary Poisson process. \\

Consider the process where customers arrive at the store, and choose products they prefer. Products can go out of stock, and become unavailable for customers who come later. The store replenish its inventory periodically. During the replenishment, the assortment can be swapped as well. The time intervals between successive replenishments may be of unequal length, but the lengths are observed. Data is recorded during the process, and is collected at replenishments. \\

Assume that the arrivals and the choice behavior in different time intervals are independent. In this case, the overall likelihood is the the product of likelihoods over all time intervals. Thus, we can model the random process of each time interval separately. Consider one time interval. Let $T$ denotes the total length of time, during which the store operates normally. Note that there is no inventory replenishment during $(0,T)$, which means assortment can only change because of stock-outs. In other words, the assortment provided to each customer is determined by the previous choices. \\

\subsection{Random Variables and Sample Paths} 

Let $N$ be the total number of customer arrivals during $[0,T]$. Let $T_1, \cdots, T_{N} \in [0,T]$ be the arrival times of the customers, and let $\mathbf{T} = (T_1, \cdots, T_{N})$. For any observed total arrivals ${n} \in \mathbb{N}_0 \defi \{0,1,\cdots \}$, let
\begin{align} 
\begin{split}
\mathcal{T}^{(1)}_{n} \defi \Big\{ \mathbf{t} \defi (t_1, \cdots, t_{n}) \ : \
		0 \leq t_1 \leq \cdots \leq t_{n} \leq T \Big\} 
\end{split} 
\end{align}
denote the set of feasible arrival time vectors. \\

Every customer makes a choice upon arrival. Let $\mathcal{A}_1$ denotes the initial product assortment that the store provides, and let $o$ denotes the "null-alternative", representing the action that a customer choose to leave the store without making any purchase. For each $a \in \mathcal{A}_1 \cup \{o\}$, let $s_a \in \mathbb{N} \defi \{1,2,\cdots\}$ be the initial stock level of product $a$. Assume that $s_a$ is finite for each $a \in \mathcal{A}_1$, which means products can go out of stock. In contrary, assume that product $o$ is always available, with "initial inventory" $s_o = \infty$. \\ 

Let $C_1, \cdots, C_N$ be the choices that customers make, and let $\mathbf{C} = (C_1, \cdots, C_{N})$. For any observed total arrivals ${n} \in \mathbb{N}_0$ during $[0,T]$, the set of feasible choice sequences is 
\begin{align} 
\begin{split}
\mathcal{C}^{(1)}_{n} \defi \Big\{ \mathbf{a} \defi (a_1, \cdots, a_{n}) \in (\mathcal{A}_1 \cup \{o\})^{n} : 
		\ \mathbb{1}(a_1 = a) + \cdots + \mathbb{1}(a_{n} = a) \leq s_a \ , \ \forall a \in \mathcal{A}_1 \Big\}
\end{split} 
\end{align}

Let $2^{\mathcal{A}_1}$ denotes the collection of subsets of $\mathcal{A}_1$. For any $i \in \{2, 3, \cdots \}$ and any $(a_1, \cdots, a_{i-1}) \in \mathcal{C}_{i-1}^{(1)}$, let $\mathcal{A}_{i} : \mathcal{C}_{i-1}^{(1)} \mapsto 2^{\mathcal{A}_1}$ be functions such that
\begin{align} 
\begin{split}
\mathcal{A}_i (a_1, \cdots a_{i-1}) &\defi \Big\{ a \in \mathcal{A}_1 : s_a > \sum_{i'=1}^{i-1} \mathbb{1}(a_{i'} = a) \Big\} \\
\end{split} 
\end{align}
Intuitively, $\mathcal{A}_i$ is the set of products left after $i-1$ choices have been made. \\

Let $(\Omega, \mathcal{F}, \mathbb{P})$ be the probability space of the random process during $[0,T]$. In this paper, we assume that $T$, $\mathcal{A}_1$ and $s_a, \ \forall a \in \mathcal{A}_1$ are known input. \\

\subsection{The Probability Model} 

$N$ follows a Poisson distribution with mean $T \lambda$. For any ${n} \in \mathbb{N}_0$, we have
\begin{align} 
\begin{split}
\mathbb{P}[N = n](\lambda) = \frac{(T \lambda)^n \cdot e^{-T \lambda}}{(n)!} 
\end{split} 
\end{align} 
Meanwhile, given $N = n$ for any ${n} \in \mathbb{N}$, we have $\mathbf{T} \in \mathcal{T}^{(1)}_{n}$. In this case, the vector of arrival times $\mathbf{T}$ has the same distribution as the order statistics of $n$ i.i.d. points uniformly distributed on $[0,T]$. Therefore, for any Lebesgue measurable set $\mathcal{T} \subseteq \mathcal{T}^{(1)}_{n}$, 
\begin{align} 
\begin{split}
\mathbb{P}[\mathbf{T} \in \mathcal{T} | N = n]
		= \mathbb{P}[\mathbf{T} \in \mathcal{T} | \mathbf{T} \in \mathcal{T}^{(1)}_{n}]
		= \int_{\mathcal{T}} \Big( \frac{(n)!}{T^n} \Big) d\mathbf{t}
		= \frac{(n)!}{T^n} \cdot \mathcal{L}(\mathcal{T})
\end{split} 
\end{align} 
where $\mathcal{L}$ is the Lebesgue measure. Let $\mathcal{L}(\mathcal{T}_0^{(1)}) = 1$ by convention. \\

Assume that the choice probability of a customer only depends on the choice and the choice set. For any $\mathcal{A} \subset \mathcal{A}_1 \cup \{o\}$ and $a \in \mathcal{A}$, let $P_{a:\mathcal{A}} (\beta)$ be the probability that a customer chooses $a$ out of $\mathcal{A}$, where $\beta$ is the parameter vector of the choice model. Given $N = n$ for any ${n} \in \mathbb{N}$, we have $\mathbf{C} \in \mathcal{C}^{(1)}_{n}$, and
\begin{align} 
\begin{split}
&\quad \mathbb{P}[\mathbf{C} = \mathbf{a} | N = n] (\beta) 
	= \mathbb{P}[\mathbf{C} = \mathbf{a} | \mathbf{C} \in \mathcal{C}^{(1)}_{n}] (\beta) \\
	&= \mathbb{P}[C_1 = a_1] (\beta) \cdot \mathbb{P}[C_2 = a_2 | C_1 = a_1] (\beta) 
			\cdots \mathbb{P}[C_n = a_n | C_{n-1} = a_{n-1}, \cdots, C_1 = a_1] (\beta) \\
	&= \prod_{i=1}^n P_{a_i:\mathcal{A}_i(a_1, \cdots, a_{i-1}) \cup \{o\}} (\beta) 
\end{split} \label{measure-seq}
\end{align}
Thus, for any set $\mathcal{C} \subseteq \mathcal{C}^{(1)}_{n}$, we have
\begin{align} 
\begin{split}
\mathbb{P}[\mathbf{C} \in \mathcal{C} | \mathbf{C} \in \mathcal{C}^{(1)}_{n}] (\beta) 
		= \sum_{\mathbf{a} \in \mathcal{C}} \Big( \prod_{i=1}^n P_{a_i:\mathcal{A}_i(a_1, \cdots, a_{i-1}) \cup \{o\}} (\beta) \Big)
\end{split} 
\end{align} \\

\subsection{Likelihood of Complete Data} 

A complete observation from the random process is one specific outcome of the process. Consider any $n \in \mathbb{N}$, any $\mathbf{a} = (a_1, \cdots, a_n) \in \mathcal{C}_{n}^{(1)}$ and any $\mathbf{t} = (t_1, \cdots, t_n) \in \mathcal{T}_{n}^{(1)}$. Let $L_1(n, \mathbf{a}, \mathbf{t})$ be the likelihood that we observe $N=n$, $\mathbf{C} = \mathbf{a}$ and $\mathbf{T} = \mathbf{t}$. To derive $L_1(n, \mathbf{a}, \mathbf{t})$, consider the event 
\begin{align} 
\begin{split}
E_1(n, \mathbf{a}, \mathcal{T}) \defi \Big\{ \omega \in \Omega \ : \ N(\omega) = n \ ; \ \mathbf{T}(\omega) \in \mathcal{T} \ ; \ \mathbf{C}(\omega) = \mathbf{a} \Big\}
\end{split} 
\end{align} 
where $\mathcal{T}$ is any Lebesgue measurable subset of $\mathcal{T}_{n}^{(1)}$. We have
\begin{align} 
\begin{split}
\mathbb{P}[E_1(n, \mathbf{a}, \mathcal{T})](\lambda, \beta) &= 
	\mathbb{P}[N = n](\lambda) \cdot \mathbb{P}[\mathbf{T} \in \mathcal{T} | N = n]\cdot \mathbb{P}[\mathbf{C} = \mathbf{a} | N = n ] (\beta) \\
&= \lambda^{n} e^{-T \lambda} \cdot \mathcal{L}(\mathcal{T}) \cdot \prod_{i=1}^{n} P_{a_i:\mathcal{A}_i(a_1, \cdots, a_{i-1}) \cup \{o\}} (\beta)
\end{split} 
\end{align}
The fact that $\mathbb{P}[E_1(n, \mathbf{a}, \mathcal{T})](\lambda, \beta) = \int_{\mathcal{T}} L_1[n, \mathbf{a}, \mathbf{t}] (\lambda, \beta) \ d\mathbf{t}$ implies
\begin{align} 
\begin{split}
L_1[n, \mathbf{a}, \mathbf{t}] (\lambda, \beta) 
&= \lambda^{n} e^{-T \lambda} \cdot \prod_{i=1}^{n} P_{a_i:\mathcal{A}_i(a_1, \cdots, a_{i-1}) \cup \{o\}} (\beta) \\
\end{split} 
\end{align} \\

\subsection{Likelihood of Complete Choice Sequence} 

Consider any $n \in \mathbb{N}$ and any $\mathbf{a} = (a_1, \cdots, a_n) \in \mathcal{C}_{n}^{(1)}$. Let 
\begin{align} 
\begin{split}
E_2(\tilde{n}, \mathbf{\tilde a}) \defi \Big\{ \omega \in \Omega \ : \ N(\omega) = n \ ; \ \mathbf{C}(\omega) = \mathbf{a} \Big\}
\end{split} 
\end{align} 
be the event that we observe $N = n$ and $\mathbf{C} = \mathbf{a}$. The likelihood of observing this event is
\begin{align} 
\begin{split}
L_2[n, \mathbf{a}] (\lambda, \beta) &= \mathbb{P}[E_2(\tilde{n}, \mathbf{\tilde a})](\lambda, \beta)
		= \mathbb{P}[N = n](\lambda) \cdot \mathbb{P}[\mathbf{C} = \mathbf{a} | N = n ] (\beta) \\
&= \frac{(T \lambda)^{n} e^{-T \lambda}}{(n)!} \cdot \prod_{i=1}^{n} P_{a_i:\mathcal{A}_i(a_1, \cdots, a_{i-1}) \cup \{o\}} (\beta) \\
\end{split} 
\end{align} 

It is clear that $L_2[n, \mathbf{a}] (\lambda, \beta) \propto L_1[n, \mathbf{a}, \mathbf{t}] (\lambda, \beta)$ with respect to $\lambda$ and $\beta$. In other words, given the total number of arrivals, when a complete choice sequence is observed, observing the value of $\mathbf{T}$ does not provide new information about $\lambda$ or $\beta$. This holds true since $\mathbf{T}$ and $\mathbf{C}$ are independent given $N$. \\

Besides, note that both $L_1[n, \mathbf{a}, \mathbf{t}] (\lambda, \beta)$ and $L_2[n, \mathbf{a}] (\lambda, \beta)$ can be separated into two parts, with one (arrival likelihood) depends only on $\lambda$, while the other (choice likelihood) only depends on $\beta$. In other words, when the complete choice sequence is given, the total arrival is irrelevant to the estimation of the choice model parameters. \\

\section{Transaction Data}

Let $\tilde{N} \defi \sum_{i=1}^{N} \mathbb{1}(C_i \neq o)$ be the total number of transactions (real purchases) that customers make during $[0,T]$. Let $\tilde{T}_1, \cdots, \tilde{T}_{\tilde{N}} \in [0,T]$ be the transaction timestamps of the customers, defined as
\begin{align} 
\begin{split}
\tilde{T}_i &\defi T_{\min\{ \tilde{n}_i \in \mathbb{N} \ : \ \mathbb{1}(C_1 \neq o) + \cdots + \mathbb{1}(C_{\tilde{n}_i} \neq o) \ = \ i \}} 
\ , \ i = 1, \cdots, \tilde{N} \\
\end{split} 
\end{align}
and let $\mathbf{\tilde T} \defi (\tilde{T}_1, \cdots, \tilde{T}_{\tilde{N}})$. For any observed total number of transactions $\tilde{n} \in \mathbb{N}_0$, the set of all feasible arrival time vectors $\mathbf{\tilde T}$ is $\mathcal{T}^{(1)}_{\tilde{n}}$. Similarly, let $\tilde{C}_1, \cdots, \tilde{C}_{\tilde{N}}$ be the purchases that customers make, defined as 
\begin{align} 
\begin{split}
\tilde{C}_i &\defi C_{\min\{ \tilde{n}_i \in \mathbb{N} \ : \ \mathbb{1}(C_1 \neq 0) + \cdots + \mathbb{1}(C_{\tilde{n}_i} \neq 0) \ = \ i \}} 
\ , \ i = 1, \cdots, \tilde{N} \\
\end{split} 
\end{align}
and let $\mathbf{\tilde C} \defi (\tilde{C}_1, \cdots, \tilde{C}_{\tilde{N}})$. For any observed total number of transactions $\tilde{n} \in \mathbb{N}_0$ during $[0,T]$, the set of feasible choice sequences is 
\begin{align} 
\begin{split}
\mathcal{\tilde C}^{(1)}_{\tilde{n}} \defi \Big\{ \mathbf{ \tilde a} \defi (\tilde{a}_1, \cdots, \tilde{a}_{\tilde{n}}) \in (\mathcal{A}_1)^{\tilde{n}} : 
		\ \mathbb{1}(\tilde{a}_1 = a) + \cdots + \mathbb{1}(\tilde{a}_{\tilde{n}} = a) \leq s_a \ , \ \forall a \in \mathcal{A}_1 \Big\}
\end{split} 
\end{align}
Let $\mathcal{\tilde A}_1 \defi \mathcal{A}_1$. For any $i \in \{2, 3, \cdots \}$ and $(\tilde{a}_1, \cdots, \tilde{a}_{i-1}) \in \mathcal{\tilde C}_{i-1}^{(1)}$, define functions $\mathcal{\tilde A}_{i} : \mathcal{\tilde C}_{i-1}^{(1)} \mapsto 2^{\mathcal{A}_1}$ as
\begin{align} 
\begin{split}
\mathcal{\tilde A}_i (\tilde{a}_1, \cdots, \tilde{a}_{i-1}) &\defi \Big\{ a \in \mathcal{A}_1 : s_a > \sum_{i'=1}^{i-1} \mathbb{1}(\tilde{a}_{i'} = a) \Big\} \\
\end{split} 
\end{align}
Intuitively, $\mathcal{\tilde A}_i$ is the set of products left after $i-1$ transactions have been made. \\

Let $N_o \defi N - \tilde{N}$ be the number of customers who choose to leave without making a purchase during $[0,T]$. For each $a \in \mathcal{A}_1$, let $N_a \defi \sum_{i=1}^{\tilde{N}} \mathbb{1}(\tilde{C}_i = a)$ be the total number of product $a$ sold during $[0,T]$. Let $K \defi \sum_{a \in \mathcal{A}_1} \mathbb{1} ( N_a = s_a ) $ denotes the number of products that are out of stock by time $T$. \\

\subsection{Additional Functions and Random Variables} 

For any $\mathbf{a} \in \cup_{n=1}^{\infty} \ \mathcal{C}^{(1)}_{n}$ and any $j \in \mathbb{N}$, let $l^{[j]}(\mathbf{a}) : \cup_{n=1}^{\infty} \ \mathcal{C}^{(1)}_{n} \mapsto \mathbb{N}$ be a function such that 
\begin{align} 
\begin{split}
l^{[j]}(\mathbf{a}) &\defi \max \ \Big\{ i \in \{1,\cdots, |\mathbf{a}|+1 \} : 
		|\mathcal{A}_i (a_1, \cdots, a_{i-1})| \geq |\mathcal{A}_1| - j + 1 \Big\} 
\end{split} 
\end{align} 
Intuitively, $l^{[j]}(\mathbf{a})$ is the index of arrivals where the $j$-th stock-out happens. For consistency, we let $l^{[0]}(\mathbf{a}) \defi 0$. \\

Similarly, for any $\mathbf{\tilde a} \in \cup_{\tilde{n}=1}^{\infty} \ \mathcal{\tilde C}^{(1)}_{\tilde{n}}$ and any $j \in \mathbb{N}$, let $\tilde{l}^{[j]}(\mathbf{\tilde a}) : \cup_{\tilde{n}=1}^{\infty} \ \mathcal{\tilde C}^{(1)}_{\tilde{n}} \mapsto \mathbb{N}$ be a function such that 
\begin{align} 
\begin{split}
\tilde{l}^{[j]}(\mathbf{\tilde a}) &\defi \max \ \Big\{ i \in \{1,\cdots, |\mathbf{\tilde a}|+1 \} : 
		|\mathcal{\tilde A}_i (\tilde{a}_1, \cdots, \tilde{a}_{i-1})| \geq |\mathcal{A}_1| - j + 1 \Big\} 
\end{split} 
\end{align} 
Intuitively, $\tilde{l}^{[j]}(\mathbf{\tilde a})$ is the index of transactions where the $j$-th stock-out happens. For consistency, we let $\tilde{l}^{[0]}(\mathbf{\tilde a}) \defi 0$. \\

Then, let $\mathcal{\check A}^{[j]} : \cup_{\tilde{n}=1}^{\infty} \ \mathcal{\tilde C}^{(1)}_{\tilde{n}} \mapsto 2^{\mathcal{A}_1}$ be the $j$-th distinct assortment provided, defined as
\begin{align} 
\begin{split}
\mathcal{\check A}^{[j]}(\mathbf{\tilde a}) &\defi 
		\mathcal{\tilde A}_{\tilde{l}^{[j-1]}(\mathbf{\tilde a})+1} (\tilde{a}_1, \cdots, \tilde{a}_{\tilde{l}^{[j-1]}(\mathbf{\tilde a})}) = \cdots =
		\mathcal{\tilde A}_{\tilde{l}^{[j]}(\mathbf{\tilde a})} (\tilde{a}_1, \cdots, \tilde{a}_{\tilde{l}^{[j]}(\mathbf{\tilde a}) -1}) 
\end{split} 
\end{align} 

Let $\tilde{T}_0 \defi 0$ and $\tilde{T}_{\tilde{N} + 1} \defi T$. For $j = 1, \cdots, K+1$, let
\begin{align} 
\begin{split}
\check{T}^{[j]} &\defi \tilde{T}_{\tilde{l}^{[j]}(\mathbf{\tilde C})} - \tilde{T}_{\tilde{l}^{[j-1]}(\mathbf{\tilde C})}
\end{split} 
\end{align} 
denotes the length of time that $\mathcal{\check A}^{[j]}(\mathbf{\tilde C})$ is provided. Let $\check{\mathbf{T}} \defi (\check{T}^{[1]}, \cdots, \check{T}^{[K+1]})$. \\


Finally, let $N^{[j]} \defi l^{[j]}(\mathbf{C}) - l^{[j-1]}(\mathbf{C}) - 1$ denotes the number of arrivals that encounter the $j$-th distinct assortment, excluding the one that leads to a stock-out (if any). Similarly, let $\tilde{N}^{[j]} \defi \tilde{l}^{[j]}(\mathbf{\tilde C}) - \tilde{l}^{[j-1]}(\mathbf{\tilde C}) - 1$, and let $N_o^{[j]} \defi N^{[j]} - \tilde{N}^{[j]}$. For convenience, we let $\mathbf{N}_o \defi (N_o^{[1]}, \cdots, N_o^{[K+1]})$. \\

\subsection{Likelihood of a Transaction Sequence with Transaction Timestamps} 

Consider any $\tilde{n} \in \mathbb{N}$, any $\mathbf{\tilde a} = (\tilde{a}_1, \cdots, \tilde{a}_{\tilde{n}}) \in \mathcal{\tilde C}_{\tilde{n}}^{(1)}$ and any $\mathbf{\tilde t} = (\tilde{t}_1, \cdots, \tilde{t}_{\tilde{n}}) \in \mathcal{T}_{\tilde{n}}^{(1)}$. Let $L_3[\tilde{n}, \mathbf{\tilde a}, \mathbf{\tilde t}] (\lambda, \beta)$ be the likelihood that we observe $\tilde{N}=\tilde{n}$, $\mathbf{\tilde C} = \mathbf{\tilde a}$ and $\mathbf{\tilde T} = \mathbf{\tilde t}$. Note that given $\mathbf{\tilde C} = \mathbf{\tilde a}$, the value of $K$ is known. Let it be $k$. Similarly, given $\mathbf{\tilde C} = \mathbf{\tilde a}$ and $\mathbf{\tilde T} = \mathbf{\tilde t}$, the values of $\check{T}^{[1]}, \cdots, T^{[k+1]}$ and $\tilde{N}^{[1]}, \cdots, \tilde{N}^{[k+1]}$ are known. Let them be $\check{t}^{[1]}, \cdots, \check{t}^{[k+1]}$ and $\tilde{n}^{[1]}, \cdots, \tilde{n}^{[k+1]}$. \\

Intuitively, the structure of $L_3[\tilde{n}, \mathbf{\tilde a}, \mathbf{\tilde t}] (\lambda, \beta)$ should be very similar to $L_1[n, \mathbf{a}, \mathbf{t}] (\lambda, \beta)$, except that we are working with transaction (true purchases) instead of arrivals. During the $j$-th time interval, the arrival rate of transactions is be $(1 - P_{o: \mathcal{\check A}^{[j]}(\mathbf{\tilde a}) \cup \{o\} } (\beta)) \cdot \lambda$. Thus, we should have
\begin{align} 
\begin{split}
L_3[\tilde{n}, \mathbf{\tilde a}, \mathbf{\tilde t}] (\lambda, \beta) 
&= \Big( \prod_{i=1}^{\tilde{n}} 
		\frac{P_{\tilde{a}_i:\mathcal{\tilde A}_{i}(\tilde{a}_1, \cdots, \tilde{a}_{i-1}) \cup \{o\}} (\beta)}
				{1 - P_{o:\mathcal{\tilde A}_{i}(\tilde{a}_1, \cdots, \tilde{a}_{i-1}) \cup \{o\}} (\beta)} \Big) 
		\cdot \prod_{j=1}^{k+1} \ 
		\Big( (1 - P_{o: \mathcal{\check A}^{[j]}(\mathbf{\tilde a}) \cup \{o\} } (\beta)) \cdot \lambda \Big)^{\tilde{n}^{[j]} + \mathbb{1}(j<k+1)} \cdot 
		e^{-(1 - P_{o: \mathcal{\check A}^{[j]}(\mathbf{\tilde a}) \cup \{o\} } (\beta)) \cdot \check{t}^{[j]} \lambda} 
\\
&= \Big( \prod_{i=1}^{\tilde{n}} 
		P_{\tilde{a}_i:\mathcal{\tilde A}_{i}(\tilde{a}_1, \cdots, \tilde{a}_{i-1}) \cup \{o\}} (\beta) \Big) \cdot 
		\lambda^{\tilde{n}} \cdot e^{-\sum_{j=1}^{k+1} (1 - P_{o: \mathcal{\check A}^{[j]}(\mathbf{\tilde a}) \cup \{o\} } (\beta)) \cdot \check{t}^{[j]} \lambda} 
\end{split} \label{mle-3}
\end{align}

To validate (\ref{mle-3}), let us consider the event 
\begin{align} 
\begin{split}
E_3(\tilde{n}, \mathbf{\tilde a}, \mathcal{\tilde T}) &\defi 
		\Big\{ \omega \in \Omega \ : \ \tilde{N}(\omega) = \tilde{n} \ , \ 
		\mathbf{\tilde T}(\omega) \in \mathcal{\tilde T} \ , \ 
		\mathbf{\tilde C}(\omega) = \mathbf{\tilde a} \Big\} 
\end{split} 
\end{align} 
where $\mathcal{\tilde T}$ is a Lebesgue measurable subset of $\mathcal{T}_{\tilde{n}}^{(1)}$. \\ 

Since stock-out changes the assortment (and thus the arrival rate of transactions), the probability measure of $\tilde{N}$ is not immediate. Meanwhile, even when the value of $\tilde{N}$ is given, $\mathbf{\tilde T}$ and $\mathbf{\tilde C}$ are still not independent. Thus, to work with $E_3(\tilde{n}, \mathbf{\tilde a}, \mathcal{\tilde T})$, we need to transform it into
\begin{align} 
\begin{split}
E_3(\tilde{n}, \mathbf{\tilde a}, \mathcal{\tilde T}) &= 
		\Big\{ \omega \in \Omega \ : \ 
		\mathbf{C}(\omega) \in \mathcal{C}_{\tilde{n}, \mathbf{\tilde a}}^{(2)} \ , \ 
		\mathbf{T}(\omega) \in \bigcup_{\mathbf{\tilde t} \in \mathcal{\tilde T}} 
		\mathcal{T}_{\tilde{n}, \mathbf{\tilde t}, \mathbf{C}(\omega)}^{(2)} \Big\}
\end{split} 
\end{align}
where
\begin{align} 
\begin{split}
\mathcal{C}_{\tilde{n}, \mathbf{\tilde a}}^{(2)} \defi \bigcup_{n = \tilde{n}}^\infty \Big\{ \mathbf{a} \in \mathcal{C}_{n}^{(1)} :  
	\quad \exists \quad & (m_1, \cdots, m_{\tilde{n}}) \in \{1, \cdots, n\}^{\tilde{n}} \\
	\text{s.t.} \quad & m_1 < \cdots < m_{\tilde{n}} \\
	\text{and} \quad & a_{m_i} = \tilde{a}_i \quad , \quad \forall i=1, \cdots, \tilde{n} \\
	\text{and} \quad & a_{i} = o \quad , \quad \forall i \in 
			\{1, \cdots, n\} \setminus \{m_1, \cdots, m_{\tilde{n}}\} \ \Big\} 
\end{split} 
\end{align}
is the set of choice sequences $\mathbf{a}$ in $\cup_{n=\tilde{n}}^\infty \ \mathcal{C}_{n}^{(1)}$ that are consistent (feasible) with $\mathbf{\tilde a}$, and 
\begin{align} 
\begin{split}
\mathcal{T}_{\tilde{n}, \mathbf{\tilde t}, \mathbf{a}}^{(2)} \defi \Big\{ \mathbf{t} \in \mathcal{T}_{|\mathbf{a}|}^{(1)} :  
	\quad \exists \quad & (m_1, \cdots, m_{\tilde{n}}) \in \{1, \cdots, |\mathbf{a}|\}^{\tilde{n}} \\
	\text{s.t.} \quad & m_1 < \cdots < m_{\tilde{n}} \\
	\text{and} \quad & t_{m_i} = \tilde{t}_i \quad , \quad \forall i=1, \cdots, \tilde{n} \\
	\text{and} \quad & a_{i} = o \quad , \quad \forall i \in 
			\{1, \cdots, |\mathbf{a}|\} \setminus \{m_1, \cdots, m_{\tilde{n}}\} \ \Big\} \\
\end{split} 
\end{align}
is the set of $\mathbf{t}$ in $\cup_{n=\tilde{n}}^\infty \ \mathcal{T}_{n}^{(1)}$ that are consistent with $\mathbf{\tilde t}$ and any $\mathbf{a} \in \mathcal{C}_{\tilde{n}, \mathbf{\tilde a}}^{(2)}$, We have
\begin{align} 
\begin{split}
&\mathbb{P}[E_3(\tilde{n}, \mathbf{\tilde a}, \mathcal{\tilde T})] (\lambda, \beta) 
= \sum_{\mathbf{a} \in \mathcal{C}_{\tilde{n}, \mathbf{\tilde a}}^{(2)}} 
		\mathbb{P}[N = |\mathbf{a}|] (\lambda) \cdot 
		\mathbb{P}[\mathbf{C} = \mathbf{a} | N = |\mathbf{a}|] (\beta) \cdot 
		\mathbb{P}[\mathbf{T} \in \bigcup_{\mathbf{\tilde t} \in \mathcal{\tilde T}} \mathcal{T}_{\tilde{n}, \mathbf{\tilde t}, \mathbf{a}}^{(2)} | N = |\mathbf{a}|] \\
\end{split} 
\end{align}

The set $\mathcal{C}_{\tilde{n}, \mathbf{\tilde a}}^{(2)}$ contains an infinite number of elements, which makes it hard to evaluate $\mathbb{P}[E_3(\tilde{n}, \mathbf{\tilde a}, \mathcal{\tilde T})] (\lambda, \beta)$. To simplify the formulation, for any $\mathbf{n}_o \defi (n_o^{[1]}, \cdots, n_o^{[k+1]})$ such that $\mathbf{n}_o \in \mathbb{N}_0^{k+1}$, let
\begin{align} 
\begin{split}
\mathcal{C}_{\tilde{n}, \mathbf{\tilde a}, \mathbf{n}_o}^{(3)} \defi \Big\{ \mathbf{a} \in \mathcal{C}_{\tilde{n} + \sum_{j=1}^{k+1} n_o^{[j]}}^{(1)} :  
	\quad \exists \quad & (m_1, \cdots, m_{\tilde{n}}) \in \{1, \cdots, \tilde{n} + \sum_{j=1}^{k+1} n_o^{[j]}\}^{\tilde{n}} \\
	\text{s.t.} \quad & m_1 < \cdots < m_{\tilde{n}} \\
	\text{and} \quad & a_{m_i} = \tilde{a}_i \quad , \quad \forall i=1, \cdots, \tilde{n} \\
	\text{and} \quad & a_{i} = o \quad , \quad \forall i \in \{1, \cdots, \tilde{n} + \sum_{j=1}^{k+1} n_o^{[j]}\} \setminus \{m_1, \cdots, m_{\tilde{n}}\} \\
	\text{and} \quad & m_{\tilde{l}^{[j]}(\mathbf{\tilde a})} 
			= \tilde{l}^{[j]}(\mathbf{\tilde a}) + \sum_{j'=1}^j n_o^{[j]} \quad , \quad j = 1, \cdots, k \ \Big\} \\
\end{split} 
\end{align}
By (\ref{measure-seq}), for any $\mathbf{a} \in \mathcal{C}_{\tilde{n}, \mathbf{\tilde a}, \mathbf{n}_o}^{(3)}$, we have
\begin{align} 
\begin{split}
\mathbb{P} \Big[\mathbf{C} = \mathbf{a} \Big| N = \tilde{n} + \sum_{j=1}^{k+1} n_o^{[j]} \Big]  
= \Big( \prod_{i=1}^{\tilde{n}} P_{\tilde{a}_i:\mathcal{\tilde A}_{i} (\tilde{a}_1, \cdots, \tilde{a}_{i-1}) \cup \{o\}} (\beta) \Big) 
		\prod_{j=1}^{k+1} \Big( P_{o: \mathcal{\check A}^{[j]}(\mathbf{\tilde a}) \cup \{o\} } (\beta) \Big)^{n_o^{[j]}} 
\end{split} 
\end{align} 
which is a constant of $\beta$. Since the $\mathcal{C}_{\tilde{n}, \mathbf{\tilde a}, \mathbf{n}_o}^{(3)}$ sets are non-overlapping, we can partition $\mathcal{C}_{\tilde{n}, \mathbf{\tilde a}}^{(2)}$ as 
\begin{align} 
\begin{split}
\mathcal{C}_{\tilde{n}, \mathbf{\tilde a}}^{(2)} = 
		\bigcup_{\mathbf{n}_o \in \mathbb{N}_0^{k+1}} 
		\mathcal{C}_{\tilde{n}, \mathbf{\tilde a}, \mathbf{n}_o}^{(3)}
\end{split} 
\end{align} 
Therefore
\begin{align} 
\begin{split}
&\quad \mathbb{P}[E_3(\tilde{n}, \mathbf{\tilde a}, \mathcal{\tilde T})] (\lambda, \beta) \\
&= \sum_{\mathbf{n}_o \in \mathbb{N}_0^{k+1}} 
		\sum_{\mathbf{a} \in \mathcal{C}_{\tilde{n}, \mathbf{\tilde a}, \mathbf{n}_o}^{(3)}} 
		\mathbb{P}\Big[N = \tilde{n} + \sum_{j=1}^{k+1} n_o^{[j]} \Big] (\lambda) \cdot 
		\mathbb{P}\Big[\mathbf{C} = \mathbf{a} \Big| N = \tilde{n} + \sum_{j=1}^{k+1} n_o^{[j]} \Big] (\beta) \cdot 
		\mathbb{P}\Big[\mathbf{T} \in \bigcup_{\mathbf{\tilde t} \in \mathcal{T}} 
			\mathcal{T}_{\tilde{n}, \mathbf{\tilde t}, \mathbf{a}}^{(2)} \Big| N = \tilde{n} + \sum_{j=1}^{k+1} n_o^{[j]} \Big] 
\\
&= \sum_{\mathbf{n}_o \in \mathbb{N}_0^{k+1}} 
		\mathbb{P}\Big[N = \tilde{n} + \sum_{j=1}^{k+1} n_o^{[j]} \Big] (\lambda) \cdot 
		\Big( \prod_{i=1}^{\tilde{n}} P_{\tilde{a}_i:\mathcal{\tilde A}_{i}(\tilde{a}_1, \cdots, \tilde{a}_{i-1}) \cup \{o\}} (\beta) \Big) 
			\prod_{j=1}^{k+1} \Big( P_{o: \mathcal{\check A}^{[j]}(\mathbf{\tilde a}) \cup \{o\} } (\beta) \Big)^{n_o^{[j]}} \cdot 
		\mathbb{P}\Big[\mathbf{T} \in \bigcup_{\mathbf{\tilde t} \in \mathcal{T}} 
			\bigcup_{\mathbf{a} \in \mathcal{C}_{\tilde{n}, \mathbf{\tilde a}, \mathbf{n}_o}^{(3)}} 
			\mathcal{T}_{\tilde{n}, \mathbf{\tilde t}, \mathbf{a}}^{(2)} | N = \tilde{n} + \sum_{j=1}^{k+1} n_o^{[j]} \Big] 
\end{split}
\end{align} 

Now, note that for any $\mathbf{n}_o \in \mathbb{N}_0^{k+1}$, 
\begin{align} 
\begin{split}
\bigcup_{\mathbf{a} \in \mathcal{C}_{\tilde{n}, \mathbf{\tilde a}, \mathbf{n}_o}^{(3)}} \mathcal{T}_{\tilde{n}, \mathbf{\tilde t}, \mathbf{a}}^{(2)}
= \Big\{ \mathbf{t} \in \mathcal{T}_{n}^{(1)} :  
	\quad \exists \quad & (m_1, \cdots, m_{\tilde{n}}) \in \{1, \cdots, \tilde{n} + \sum_{j=1}^{k+1} n_o^{[j]}\}^{\tilde{n}} \\
	\text{s.t.} \quad & m_1 < \cdots < m_{\tilde{n}} \\
	\text{and} \quad & t_{m_i} = \tilde{t}_i \quad , \quad \forall i=1, \cdots, \tilde{n} \\
	\text{and} \quad & m_{\tilde{l}^{[j]}(\mathbf{\tilde a})} = \tilde{l}^{[j]}(\mathbf{\tilde a}) + \sum_{j'=1}^{j} n_o^{[j']} \quad , \quad j = 1, \cdots, k \ \Big\} \\
\end{split} \label{set-r3}
\end{align} 
In other words, given $\tilde{N} = \tilde{n}$, $\mathbf{\tilde T} = \mathbf{\tilde t}$, $\mathbf{\tilde C} = \mathbf{\tilde c}$ and $\mathbf{N}_o = \mathbf{n}_o$, for each $j=1,\cdots, k+1$, there are $n_o^{[j]}$ arrivals associated with the null alternative $o$, that has the same distribution as the order statistics of $n_o^{[j]}$ i.i.d. points uniformly distributed on the $j$-th time interval with length $\check{t}^{[j]}$. It implies that
\begin{align} 
\begin{split}
\mathbb{P}\Big[\mathbf{T} \in \bigcup_{\mathbf{\tilde t} \in \mathcal{T}} 
		\bigcup_{\mathbf{a} \in \mathcal{C}_{\tilde{n}, \mathbf{\tilde a}, \mathbf{n}_o}^{(3)}} 
		\mathcal{T}_{\tilde{n}, \mathbf{\tilde t}, \mathbf{a}}^{(2)} \Big| N = \tilde{n} + \sum_{j=1}^{k+1} n_o^{[j]} \Big] 
&= \int_{\mathcal{T}} \Big( \frac{(\tilde{n} + \sum_{j=1}^{k+1} n_o^{[j]})!}{T^{\tilde{n} + \sum_{j=1}^{k+1} n_o^{[j]} }} \cdot 
		\prod_{j=1}^{k+1} \frac{(\check{t}^{[j]})^{n_o^{[j]}}}{(n_o^{[j]})!} \Big) \ d\mathbf{\tilde t}
\end{split} 
\end{align} 
Therefore
\begin{align} 
\begin{split}
&\quad \mathbb{P}[E_3(\tilde{n}, \mathbf{\tilde a}, \mathcal{\tilde T})] (\lambda, \beta) \\
\\
&= \sum_{\mathbf{n}_o \in \mathbb{N}_0^{k+1}} 
		\Big( \frac{(T \lambda)^{\tilde{n} + \sum_{j=1}^{k+1} n_o^{[j]}} \cdot e^{-T \lambda}}{({\tilde{n} + \sum_{j=1}^{k+1} n_o^{[j]}})!} \Big) \cdot 
		\Big( \prod_{i=1}^{\tilde{n}} P_{\tilde{a}_i:\mathcal{\tilde A}_{i}(\tilde{a}_1, \cdots, \tilde{a}_{i-1}) \cup \{o\}} (\beta) \Big) 
			\prod_{j=1}^{k+1} \Big( P_{o: \mathcal{\check A}^{[j]}(\mathbf{\tilde a}) \cup \{o\} } (\beta) \Big)^{n_o^{[j]}} \cdot 
		\mathbb{P}\Big[\mathbf{T} \in \bigcup_{\mathbf{\tilde t} \in \mathcal{T}} 
			\bigcup_{\mathbf{a} \in \mathcal{C}_{\tilde{n}, \mathbf{\tilde a}, \mathbf{n}_o}^{(3)}} 
			\mathcal{T}_{\tilde{n}, \mathbf{\tilde t}, \mathbf{a}}^{(2)} | N = \tilde{n} + \sum_{j=1}^{k+1} n_o^{[j]} \Big] 
\\
&= \sum_{\mathbf{n}_o \in \mathbb{N}_0^{k+1}} \int_{\mathcal{T}} \Bigg( 
			\lambda^{\tilde{n} + \sum_{j=1}^{k+1} n_o^{[j]}} \cdot e^{-T \lambda} \cdot 
			\Big( \prod_{j=1}^{k+1} \frac{(\check{t}^{[j]})^{n_o^{[j]}}}{(n_o^{[j]})!} \Big) \cdot
			\Big( \prod_{i=1}^{\tilde{n}} P_{\tilde{a}_i:\mathcal{\tilde A}_{i}(\tilde{a}_1, \cdots, \tilde{a}_{i-1}) \cup \{o\}} (\beta) \Big) 
				\prod_{j=1}^{k+1} \Big( P_{o: \mathcal{\check A}^{[j]}(\mathbf{\tilde a}) \cup \{o\} } (\beta) \Big)^{n_o^{[j]}} \Bigg) d\mathbf{\tilde t}
\\
&= \sum_{\mathbf{n}_o \in \mathbb{N}_0^{k+1}} \int_{\mathcal{T}} \Bigg( 
			\lambda^{\tilde{n}} e^{-T \lambda} \cdot 
			\Big( \prod_{i=1}^{\tilde{n}} P_{\tilde{a}_i:\mathcal{\tilde A}_{i}(\tilde{a}_1, \cdots, \tilde{a}_{i-1}) \cup \{o\}} (\beta) \Big) 
			\Big( \prod_{j=1}^{k+1} \frac{\Big( P_{o: \mathcal{\check A}^{[j]}(\mathbf{\tilde a}) \cup \{o\} } (\beta) \cdot
				\check{t}^{[j]} \lambda \Big)^{n_o^{[j]}}}{(n_o^{[j]})!} \Big) \Bigg) d\mathbf{\tilde t}
\\
&= \int_{\mathcal{T}} \Bigg( \lambda^{\tilde{n}} e^{-T \lambda} \cdot 
			\Big( \prod_{i=1}^{\tilde{n}} P_{\tilde{a}_i:\mathcal{\tilde A}_{i}(\tilde{a}_1, \cdots, \tilde{a}_{i-1}) \cup \{o\}} (\beta) \Big) 
			\prod_{j=1}^{k+1} \sum_{\mathbf{n}_o \in \mathbb{N}_0^{k+1}} 
				\frac{\Big( P_{o: \mathcal{\check A}^{[j]}(\mathbf{\tilde a}) \cup \{o\} } (\beta) \cdot
				\check{t}^{[j]} \lambda \Big)^{n_o^{[j]}}}{(n_o^{[j]})!} \Bigg) d\mathbf{\tilde t}
\\
&= \int_{\mathcal{T}} \Bigg( \lambda^{\tilde{n}} e^{-T \lambda} \cdot 
			\Big( \prod_{i=1}^{\tilde{n}} P_{\tilde{a}_i:\mathcal{\tilde A}_{i}(\tilde{a}_1, \cdots, \tilde{a}_{i-1}) \cup \{o\}} (\beta) \Big) 
			\prod_{j=1}^{k+1} e^{ P_{o: \mathcal{\check A}^{[j]}(\mathbf{\tilde a}) \cup \{o\} } (\beta) \cdot
				\check{t}^{[j]} \lambda } \Bigg) d\mathbf{\tilde t}
\end{split}
\end{align} 

Note that $T = \sum_{j=1}^{k+1} \check{t}^{[j]}$ always holds. Therefore, 
\begin{align} 
\begin{split}
\mathbb{P}[E_3(\tilde{n}, \mathbf{\tilde a}, \mathcal{\tilde T})] (\lambda, \beta) 
&= \Big( \prod_{i=1}^{\tilde{n}} P_{\tilde{a}_i:\mathcal{\tilde A}_{i}(\tilde{a}_1, \cdots, \tilde{a}_{i-1}) \cup \{o\}} (\beta) \Big) 
			\int_{\mathcal{\tilde T}} \lambda^{\tilde{n}} \cdot 
			e^{-\sum_{j=1}^{k+1} (1 - P_{o: \mathcal{\check A}^{[j]}(\mathbf{\tilde a}) \cup \{o\} } (\beta)) \cdot \check{t}^{[j]} \lambda} 
			\ d\mathbf{\tilde t}
\end{split}
\end{align} 
Meanwhile, 
\begin{align} 
\begin{split}
\int_{\mathcal{T}} L_3[\tilde{n}, \mathbf{\tilde a}, \mathbf{\tilde t}] (\lambda, \beta) \ d\mathbf{\tilde t} = \mathbb{P}[E_3(\tilde{n}, \mathbf{\tilde a}, \mathcal{\tilde T})] (\lambda, \beta)
\end{split} 
\end{align} 
Thus, 
\begin{align} 
\begin{split}
L_3[\tilde{n}, \mathbf{\tilde a}, \mathbf{\tilde t}] (\lambda, \beta) 
&= \Big( \prod_{i=1}^{\tilde{n}} P_{\tilde{a}_i:\mathcal{\tilde A}_{i}(\tilde{a}_1, \cdots, \tilde{a}_{i-1}) \cup \{o\}} (\beta) \Big) \cdot 
			\lambda^{\tilde{n}} \cdot e^{-\sum_{j=1}^{k+1} (1 - P_{o: \mathcal{\check A}^{[j]}(\mathbf{\tilde a}) \cup \{o\} } (\beta)) \cdot \check{t}^{[j]} \lambda} 
\end{split} 
\end{align} \\

\subsection{Likelihood of a Transaction Sequence} 

Suppose that we observe $\tilde{N} = \tilde{n}$ and $\mathbf{\tilde C} = \mathbf{\tilde a}$ for some $\tilde{n} \in \mathbb{N}$ and $\mathbf{\tilde a} = (\tilde{a}_1, \cdots, \tilde{a}_{\tilde{n}}) \in \mathcal{C}_{\tilde{n}}^{(1)}$. Let 
\begin{align} 
\begin{split}
E_4(\tilde{n}, \mathbf{\tilde a}) \defi \Big\{ \omega \in \Omega \ : \ \tilde{N}(\omega) = \tilde{n} \ ; \ \mathbf{\tilde C}(\omega) = \mathbf{\tilde a} \Big\}
\end{split} 
\end{align} 
be the event that we observe $\tilde{N} = \tilde{n}$ and $\mathbf{\tilde C} = \mathbf{\tilde a}$. Again, given $\mathbf{\tilde C} = \mathbf{\tilde a}$, the values of $K$ and $\tilde{N}^{[1]}, \cdots, \tilde{N}^{[k+1]}$ are known. Let them be $k$ and $\tilde{n}^{[1]}, \cdots, \tilde{n}^{[k+1]}$. \\

\subsubsection{The Summation Representation} 

Similar to the previous part, we have
\begin{align} 
\begin{split}
E_4(\tilde{n}, \mathbf{\tilde a}) = \Big\{ \omega \in \Omega : \mathbf{C}(\omega) \in \mathcal{C}_{\tilde{n}, \mathbf{\tilde a}}^{(2)} \Big\}
\end{split} 
\end{align} 
Thus, the likelihood of observing $\tilde{N} = \tilde{n}$ and $\mathbf{\tilde C} = \mathbf{\tilde a}$ is
\begin{align} 
\begin{split}
L_4[\tilde{n}, \mathbf{\tilde a}] (\lambda, \beta) 
= \mathbb{P}[E_4(\tilde{n}, \mathbf{\tilde a})] (\lambda, \beta)
= \sum_{\mathbf{a} \in \mathcal{C}_{\tilde{n}, \mathbf{\tilde a}}^{(2)}} 
		\mathbb{P}[N = |\mathbf{a}|] (\lambda) \cdot 
		\mathbb{P}[\mathbf{C} = \mathbf{a} | N = |\mathbf{a}|] (\beta) 
\end{split} 
\end{align}
Recall that
\begin{align} 
\begin{split}
\mathcal{C}_{\tilde{n}, \mathbf{\tilde a}}^{(2)} = \bigcup_{\mathbf{n}_o \in \mathbb{N}_0^{k+1}} \mathcal{C}_{\tilde{n}, \mathbf{\tilde a}, \mathbf{n}_o}^{(3)}
\end{split} 
\end{align} 
while for any $\mathbf{a} \in \mathcal{C}_{\tilde{n}, \mathbf{\tilde a}, \mathbf{n}_o}^{(3)}$, we have
\begin{align} 
\begin{split}
\mathbb{P} \Big[\mathbf{C} = \mathbf{a} \Big| N = \tilde{n} + \sum_{j=1}^{k+1} n_o^{[j]} \Big] (\beta) 
= \Big( \prod_{i=1}^{\tilde{n}} P_{\tilde{a}_i:\mathcal{\tilde A}_{i} (\tilde{a}_1, \cdots, \tilde{a}_{i-1}) \cup \{o\}} (\beta) \Big) 
		\prod_{j=1}^{k+1} \Big( P_{o: \mathcal{\check A}^{[j]}(\mathbf{\tilde a}) \cup \{o\} } (\beta) \Big)^{n_o^{[j]}} 
\end{split} 
\end{align} 
which is a constant of $\beta$. It is clear that
\begin{align} 
\begin{split}
|\mathcal{C}_{\tilde{n}, \mathbf{\tilde a}, \mathbf{n}_o}^{(3)}| 
= \prod_{j=1}^{k+1} \frac{(n_o^{[j]} + \tilde{n}^{[j]})!}{(n_o^{[j]})! (\tilde{n}^{[j]})!}
\end{split} 
\end{align}
Therefore, we have
\begin{align} 
\begin{split}
L_4[\tilde{n}, \mathbf{\tilde a}] (\lambda, \beta) 
&= \sum_{\mathbf{n}_o \in \mathbb{N}_0^{k+1}} 
		\mathbb{P} \Big[N = \tilde{n} + \sum_{j=1}^{k+1} n_o^{[j]} \Big] (\lambda)  \cdot 
		\mathbb{P} \Big[\mathbf{C} \in \mathcal{C}_{\tilde{n}, \mathbf{\tilde a}, \mathbf{n}_o}^{(3)} \Big| N = \tilde{n} + \sum_{j=1}^{k+1} n_o^{[j]} \Big] (\beta) \\
&= \sum_{\mathbf{n}_o \in \mathbb{N}_0^{k+1}} 
		\frac{(T \lambda)^{\tilde{n} + \sum_{j=1}^{k+1} n_o^{[j]}} \cdot e^{-T \lambda}}{(\tilde{n} + \sum_{j=1}^{k+1} n_o^{[j]})!} \cdot 
		\Big( \prod_{j=1}^{k+1} \frac{(n_o^{[j]} + \tilde{n}^{[j]})!}{(n_o^{[j]})! (\tilde{n}^{[j]})!} \Big) \cdot 
		\Big( \prod_{i=1}^{\tilde{n}} P_{\tilde{a}_i:\mathcal{\tilde A}_{i} (\tilde{a}_1, \cdots, \tilde{a}_{i-1}) \cup \{o\}} (\beta) \Big) 
		\prod_{j=1}^{k+1} \Big( P_{o: \mathcal{\check A}^{[j]}(\mathbf{\tilde a}) \cup \{o\} } (\beta) \Big)^{n_o^{[j]}} 
\end{split} \label{mle-4}
\end{align} \\

\subsubsection{The Integral Representation}

First, note that $E_4(\tilde{n}, \mathbf{\tilde a}) = E_3(\tilde{n}, \mathbf{\tilde a}, \mathcal{T}_{\tilde{n}}^{(1)})$. Recall
\begin{align} 
\begin{split}
\mathbb{P}[E_3(\tilde{n}, \mathbf{\tilde a}, \mathcal{\tilde T})] (\lambda, \beta) 
&= \Big( \prod_{i=1}^{\tilde{n}} P_{\tilde{a}_i:\mathcal{\tilde A}_{i}(\tilde{a}_1, \cdots, \tilde{a}_{i-1}) \cup \{o\}} (\beta) \Big) 
			\int_{\mathcal{\tilde T}} \lambda^{\tilde{n}} \cdot 
			e^{-\sum_{j=1}^{k+1} (1 - P_{o: \mathcal{\check A}^{[j]}(\mathbf{\tilde a}) \cup \{o\} } (\beta)) \cdot \check{t}^{[j]} \lambda} 
			\ d\mathbf{\tilde t}
\end{split}
\end{align} 
For any $k' \in \mathbb{N}$, let
\begin{align} 
\begin{split}
\mathcal{\check T}_{k'}^{(1)} 
		\defi \Big\{ \mathbf{\check t} = (\check{t}^{[1]}, \cdots, \check{t}^{[k'+1]}) \ : \ \check{t}^{[1]} + \cdots + \check{t}^{[k'+1]} =T  \ ; \ \check{t}^{[j]} \geq 0 \ , \ \forall j = 1,\cdots,k'+1 \Big\}
\end{split}
\end{align} 
For any $\mathbf{\check t} \in \mathcal{\check T}_{k}^{(1)}$, let 
\begin{align} 
\begin{split}
\mathcal{T}_{\tilde{n}, \mathbf{\check t}, \mathbf{\tilde a}}^{(2)}  
		\defi \Big\{ \mathbf{\tilde t} \in \mathcal{T}_{\tilde{n}}^{(1)} \ : \ \tilde{t}_{\tilde{l}^{[j]}(\mathbf{\tilde a})} = \sum_{j'=1}^j \check{t}^{[j']} \ , \ j = 1,\cdots,k \Big\}
\end{split}
\end{align} 
We have
\begin{align} 
\begin{split}
L_4[\tilde{n}, \mathbf{\tilde a}] (\lambda, \beta) 
&= \mathbb{P}[E_3(\tilde{n}, \mathbf{\tilde a}, \mathcal{T}_{\tilde{n}}^{(1)})] (\lambda, \beta) 
\\
&= \Big( \prod_{i=1}^{\tilde{n}} P_{\tilde{a}_i:\mathcal{\tilde A}_{i}(\tilde{a}_1, \cdots, \tilde{a}_{i-1}) \cup \{o\}} (\beta) \Big) 
			\int_{\mathbf{\tilde t} \in \mathcal{T}_{\tilde{n}}^{(1)}} \lambda^{\tilde{n}} \cdot 
			e^{-\sum_{j=1}^{k+1} (1 - P_{o: \mathcal{\check A}^{[j]}(\mathbf{\tilde a}) \cup \{o\} } (\beta)) \cdot \check{t}^{[j]} \lambda} 
			\ d\mathbf{\tilde t} 
\\
&= \Big( \prod_{i=1}^{\tilde{n}} P_{\tilde{a}_i:\mathcal{\tilde A}_{i}(\tilde{a}_1, \cdots, \tilde{a}_{i-1}) \cup \{o\}} (\beta) \Big) 
			\int_{\mathbf{\tilde t} \in \bigcup_{\mathbf{\check t} \in \mathcal{\check T}_{k}^{(1)}} \mathcal{T}_{\tilde{n}, \mathbf{\check t}, \mathbf{\tilde a}}^{(2)} } \lambda^{\tilde{n}} \cdot 
			e^{-\sum_{j=1}^{k+1} (1 - P_{o: \mathcal{\check A}^{[j]}(\mathbf{\tilde a}) \cup \{o\} } (\beta)) \cdot \check{t}^{[j]} \lambda} 
			\ d\mathbf{\tilde t} 
\\
&= \Big( \prod_{i=1}^{\tilde{n}} P_{\tilde{a}_i:\mathcal{\tilde A}_{i}(\tilde{a}_1, \cdots, \tilde{a}_{i-1}) \cup \{o\}} (\beta) \Big) 
		\int_{\mathbf{\check t} \in \mathcal{\check T}_{k}^{(1)} } \lambda^{\tilde{n}} \cdot 
		e^{- \sum_{j=1}^{k+1} \big(1-P_{o: \mathcal{\check A}^{[j]}(\mathbf{\tilde a}) \cup \{o\}} (\beta)\big) \cdot \check{t}^{[j]} \lambda } 
		\cdot \mathcal{L}(\mathcal{T}_{\tilde{n}, \mathbf{\check t}, \mathbf{\tilde a}}^{(2)}) \ d \check{\mathbf{t}} 
\end{split} 
\end{align} 

The set $\mathcal{T}_{\tilde{n}, \mathbf{\check t}, \mathbf{\tilde a}}^{(2)}$ is $\tilde{n} - k$ dimensional. To evaluate its Lebesgue measure, we define
\begin{align} 
\begin{split}
\mathcal{T}_{\tilde{n}, \mathbf{\check t}, \mathbf{\tilde a}}^{(2, j)} 
		\defi \Big\{ (\tilde{t}_{\tilde{l}^{[j-1]}(\mathbf{\tilde a}) + 1}, \cdots, \tilde{t}_{\tilde{l}^{[j]}(\mathbf{\tilde a}) - 1}) \ : \ 
			\sum_{j'=1}^{j-1} \check{t}^{[j']} < \tilde{t}_{\tilde{l}^{[j-1]}(\mathbf{\tilde a}) + 1} < \cdots 
			< \tilde{t}_{\tilde{l}^{[j]}(\mathbf{\tilde a}) - 1} < \sum_{j'=1}^j \check{t}^{[j']} \Big\}
\end{split}
\end{align} 
for every $j = 1, \cdots, k+1$. It is clear that
\begin{align} 
\begin{split}
\mathcal{L}(\mathcal{T}_{\tilde{n}, \mathbf{\check t}, \mathbf{\tilde a}}^{(2)}) 
= \prod_{j=1}^{k+1} \mathcal{L}(\mathcal{T}_{\tilde{n}, \mathbf{\check t}, \mathbf{\tilde a}}^{(2, j)})
= \prod_{j=1}^{k+1} \Big( \frac{(\check{t}^{[j]})^{\tilde{n}^{[j]}}}{(\tilde{n}^{[j]})!} \Big)
= \frac{\prod_{j=1}^{k+1} (\check{t}^{[j]})^{\tilde{n}^{[j]}}}{\prod_{j=1}^{k+1} (\tilde{n}^{[j]})!} 
\end{split}
\end{align} 
which implies
\begin{align} 
\begin{split}
L_4[\tilde{n}, \mathbf{\tilde a}] (\lambda, \beta) 
&= \Big( \prod_{i=1}^{\tilde{n}} P_{\tilde{a}_i:\mathcal{\tilde A}_{i}(\tilde{a}_1, \cdots, \tilde{a}_{i-1}) \cup \{o\}} (\beta) \Big) 
		\int_{\mathbf{\check t} \in \mathcal{\check T}_{k}^{(1)} } \lambda^{\tilde{n}} \cdot 
		e^{- \sum_{j=1}^{k+1} \big(1-P_{o: \mathcal{\check A}^{[j]}(\mathbf{\tilde a}) \cup \{o\}} (\beta)\big) \cdot \check{t}^{[j]} \lambda } 
		\cdot \Big( \frac{\prod_{j=1}^{k+1} (\check{t}^{[j]})^{\tilde{n}^{[j]}}}{\prod_{j=1}^{k+1} (\tilde{n}^{[j]})!} \Big) d \check{\mathbf{t}} 
\end{split} \label{mle-4-2}
\end{align} \\

\subsubsection{An Alternative Proof of Equivalence}

The integral in (\ref{mle-4-2}) can be written as a moment generating function of a multivariate beta distribution. To see this fact, for any $(\check{t}^{[1]}, \cdots, \check{t}^{[k+1]}) \in \mathcal{\check T}_{k}^{(1)}$, let $q^{[j]} = \frac{\check{t}^{[j]}}{T}$ for $j=1, \cdots, k+1$, and let $\check{\mathbf{q}} \defi (q^{[1]}, \cdots, q^{[k+1]})$. Thus, $\check{\mathbf{q}} \in \mathbb{T}^{k}$, where $\mathbb{T}^{k} \defi \{\check{\mathbf{q}} = (q^{[1]}, \cdots, q^{[k+1]}) \geq 0 \ : \ \sum_{j=1}^{k+1} q^{[j]} = 1\}$ is the $k$ dimensional standard simplex. Now, consider the probability density function $g$ of a multivariate beta distribution (Dirichlet distribution) $MBD(\tilde{n}^{[1]} + 1, \cdots, \tilde{n}^{[k+1]} + 1)$, given by
\begin{align} 
\begin{split}
g(\check{\mathbf{q}}) = \frac{(\tilde{n}!) \prod_{j=1}^{k+1} (q^{[j]})^{\tilde{n}^{[j]}}}{\prod_{j=1}^{k+1} (\tilde{n}^{[j]})!}
\end{split} 
\end{align} 
We have
\begin{align} 
\begin{split}
L_4[\tilde{n}, \mathbf{\tilde a}] (\lambda, \beta) 
&= \Big( \prod_{i=1}^{\tilde{n}} P_{\tilde{a}_i:\mathcal{\tilde A}_{i}(\tilde{a}_1, \cdots, \tilde{a}_{i-1}) \cup \{o\}} (\beta) \Big) \cdot \frac{(T \lambda)^{\tilde{n}}}{(\tilde{n}!)} 
		\int_{\mathbb{T}^{k}} \Bigg[ e^{-T \lambda \sum_{j=1}^{k+1} \big(1-P_{o: \mathcal{\check A}^{[j]}(\mathbf{\tilde a}) \cup \{o\}} (\beta)\big) \cdot q^{[j]}} 
		\cdot \Big( \frac{(\tilde{n}!) \prod_{j=1}^{k+1} (q^{[j]})^{\tilde{n}^{[j]}}}{\prod_{j=1}^{k+1} (\tilde{n}^{[j]})!} \Big) \Bigg] d \check{\mathbf{q}} 
\\
&= \Big( \prod_{i=1}^{\tilde{n}} P_{\tilde{a}_i:\mathcal{\tilde A}_{i}(\tilde{a}_1, \cdots, \tilde{a}_{i-1}) \cup \{o\}} (\beta) \Big) \cdot 
		\frac{(T \lambda)^{\tilde{n}} \cdot e^{-T \lambda}}{(\tilde{n}!)} \cdot \mathbb{E}_{\check{\mathbf{q}} \sim g} 
		\Bigg[ e^{\ \sum_{j=1}^{k+1} \big(T \lambda \cdot P_{o: \mathcal{\check A}^{[j]}(\mathbf{\tilde a}) \cup \{o\}} (\beta)\big) \cdot q^{[j]}} \Bigg]
\end{split} \label{mle-4-3}
\end{align} 
This moment generating function can be written as the confluent form of a Lauricella hypergeometric series, 
\begin{align} 
\begin{split}
\mathbb{E}_{\check{\mathbf{q}} \sim g} 
	\Bigg[ e^{\ \sum_{j=1}^{k+1} \big(T \lambda \cdot P_{o: \mathcal{\check A}^{[j]}(\mathbf{\tilde a}) \cup \{o\}} (\beta)\big) \cdot q^{[j]}} \Bigg]
&= \sum_{\mathbf{n}_o \in \mathbb{N}_0^{k+1}} \frac{(\tilde{n})!}{(\tilde{n} + \sum_{j=1}^{k+1} n_o^{[j]})!} \cdot 
	\Big( \prod_{j=1}^{k+1} \frac{(n_o^{[j]} + \tilde{n}^{[j]})!}{(n_o^{[j]})! (\tilde{n}^{[j]})!} \Big) \cdot
	\prod_{j=1}^{k+1} \Big( T \lambda \cdot P_{o: \mathcal{\check A}^{[j]}(\mathbf{\tilde a}) \cup \{o\}} (\beta) \Big)^{n_o^{[j]}}
\end{split} \label{lauricella}
\end{align} 
By plugging (\ref{lauricella}) back into (\ref{mle-4-3}), we can get (\ref{mle-4}), which provides an alternative way to prove that (\ref{mle-4}) and (\ref{mle-4-2}) are equivalent. \\

\section{Sales Data}

For $j = 1, \cdots K$, let $\check{C}^{[j]} \defi \tilde{C}_{\tilde{l}^{[j]} (\tilde{C})}$ be the $j$-th product that is out of stock, and let $\mathbf{\check C} \defi (\check{C}^{[1]}, \cdots, \check{C}^{[K]})$. For any observed total number of stock-out events $k \in \mathbb{N}$ during $[0,T]$, the set of all feasible values of $\mathbf{\check C}$ is 
\begin{align} 
\begin{split}
\check{\mathcal{C}}_{k}^{(1)} \defi \Big\{ \check{\mathbf{a}} = (\check{a}^{[1]}, \cdots, \check{a}^{[k]}) \in \mathcal{A}^k \ : \ 
		\check{a}^{[j]} \neq \check{a}^{[j']} \ , \ \forall j, j' = 1, \cdots, k \ , \ j \neq j'  \ \Big\}
\end{split} 
\end{align} 

For each $a \in \mathcal{A}_1$ and $j = 1, \cdots, K+1$, let
\begin{align} 
\begin{split}
N_a^{[j]} &= \sum_{i=\tilde{l}^{[j-1]}(\mathbf{\tilde C})+1}^{\tilde{l}^{[j]}(\mathbf{\tilde C}) - 1} \mathbb{1}(\tilde{C}_i = a)
\end{split} 
\end{align} 
denotes the number of customers that encounter the $j$-th distinct assortment and choose $a$ out of it, excluding the one that leads to a stock-out (if any). For convenience, let $\mathbf{N}_a \defi (N_a^{[1]}, \cdots, N_a^{[K+1]})$ for all $a \in \mathcal{A}_1$, and let $\mathbf{N} \defi (N^{[1]}, \cdots, N^{[K+1]})$. \\ 

Finally, let $\mathbf{Z} \defi (N_a : a \in \mathcal{A}_1)$, and let $\vec{\mathbf{Z}} \defi (\mathbf{N}_a : a \in \mathcal{A}_1)$. \\

\subsection{Likelihood of Sales} 

Consider any $\tilde{n} \in \mathbb{N}$ and any $n_a \in \{0, 1, \cdots, s_a\}$ for each $a \in \mathcal{A}_1$, such that $\sum_{a \in \mathcal{A}_1} n_a = \tilde{n}$. Let $\mathbf{z} = (n_a : a \in \mathcal{A}_1)$, and let 
\begin{align} 
\begin{split}
E_5(\tilde{n}, \mathbf{z}) \defi \Big\{ \omega \in \Omega \ : \ N_a(\omega) = n_a \ , \ \forall a \in \mathcal{A}_1 \Big\}
\end{split} 
\end{align} 
be the event that we observe $\tilde{N} = \tilde{n}$ and $\mathbf{Z} = \mathbf{z}$. Given $\mathbf{Z} = \mathbf{z}$, the values of $K$ is known. Let it be $k$. We have
\begin{align} 
\begin{split}
E_5(\tilde{n}, \mathbf{z}) = \Big\{ \omega \in \Omega \ : \ \mathbf{C}(\omega) \in \mathcal{C}_{\tilde{n}, \mathbf{z}}^{(4)} \Big\}
\end{split} 
\end{align} 
where
\begin{align} 
\begin{split}
\mathcal{C}_{\tilde{n}, \mathbf{z}}^{(4)} \defi \Big\{ \mathbf{a} \in \bigcup_{n = \tilde{n}}^\infty \mathcal{C}_{n}^{(1)} :  
	\quad \sum_{i=1}^{|\mathbf{a}|} \mathbb{1}(a_i = a') = n_{a'} \ , \ \forall a' \in \mathcal{A}_1 \ \Big\} 
\end{split} 
\end{align}
is the set of choice sequences $\mathbf{a}$ that are consistent with $\mathbf{z}$. Thus, the likelihood of observing $E_5(\tilde{n}, \mathbf{z})$ is
\begin{align} 
\begin{split}
L_5[\tilde{n}, \mathbf{z}] (\lambda, \beta) 
&= \mathbb{P}[E_5(\tilde{n}, \mathbf{z})] (\lambda, \beta) 
= \sum_{\mathbf{a} \in \mathcal{C}_{\tilde{n}, \mathbf{z}}^{(4)}} 
		\mathbb{P}[N = |\mathbf{a}|] (\lambda) \cdot 
		\mathbb{P}[\mathbf{C} = \mathbf{a} | N = |\mathbf{a}|] (\beta) 
\end{split} 
\end{align} 

Again, to make the evaluation of $L_5[\tilde{n}, \mathbf{z}] (\lambda, \beta)$ easier, we can partition $\mathcal{C}_{\tilde{n}, \mathbf{z}}^{(4)}$ into sets, within which the likelihood of observing the choice sequence is a constant of $\beta$. To achieve that, first let 
\begin{align} 
\begin{split}
\check{\mathcal{C}}_{\mathbf{z}}^{(2)} \defi \Big\{ \check{\mathbf{a}} = (\check{a}^{[1]}, \cdots, \check{a}^{[k]}) \in \check{\mathcal{C}}_{k}^{(1)} \ : \ 
		n_{\check{a}^{[j]}} = s_{\check{a}^{[j]}} \ , \ \forall j = 1, \cdots, k \ \Big\}
\end{split} 
\end{align} 
denotes the set of all feasible values of $\check{\mathbf{C}}$ given $\mathbf{Z} = \mathbf{z}$. For any $\mathbf{\check a} \in \check{\mathcal{C}}_{\mathbf{z}}^{(2)}$, let
\begin{align} 
\begin{split}
\vec{\mathcal{Z}}_{\mathbf{\check a}, \mathbf{z}}^{(1)} \defi \Bigg\{ \vec{\mathbf{z}} \defi 
		\Big( (n_a^{[1]}, \cdots, n_a^{[k+1]}) : a \in \mathcal{A}_1 \Big)
		\in \mathbb{N}^{|\mathcal{A}_1| \times (k+1)} \ : \
		& n_{\check{a}^{[j]}}^{[j']} = 0 \ , \ \forall j, j' \in \{1,\cdots,k\} \ , \ j' > j \\
		& \sum_{j=1}^{k+1} n_{a}^{[j]} = n_a - \mathbb{1}(n_a = s_a) \ , \ \forall a \in \mathcal{A}_1 \ 
\Bigg\}
\end{split} 
\end{align} 
denotes the feasible set of value that $\vec{\mathbf{Z}}$ can take given $\mathbf{Z} = \mathbf{z}$ and $\check{\mathbf{C}} = \check{\mathbf{a}}$. Now, for any $\mathbf{\check a} \in \check{\mathcal{C}}_{\mathbf{z}}^{(2)}$, any $\vec{\mathbf{z}} \in \vec{\mathcal{Z}}_{\mathbf{z}, \mathbf{\check a}}^{(1)}$ and any $\mathbf{n}_o = (n_o^{[1]}, \cdots, n_o^{[k+1]}) \in \mathbb{N}_0^{k+1}$, let
\begin{align} 
\begin{split}
\mathcal{C}_{\tilde{n}, \mathbf{\check a}, \vec{\mathbf{z}}, \mathbf{n}_o}^{(5)} \defi \Bigg\{ 
		\mathbf{a} \in \mathcal{C}^{(1)}_{\tilde{n} + \sum_{j=1}^{k+1} n_o^{[j]}} \ : \ 
		\check{a}^{[j]} &= a_{l^{[j]}(\mathbf{a})} \ , \ j = 1, \cdots, k \ ;  \\
		n_a^{[j]} &= \sum_{i=l^{[j-1]}(\mathbf{a})+1}^{l^{[j]}(\mathbf{a}) - 1} \mathbb{1}(a_i = a) \ , \ j = 1, \cdots, k+1 \ , \ \forall a \in \mathcal{A}_1 \ ;  \\
		n_o^{[j]} &= \sum_{i=l^{[j-1]}(\mathbf{a})+1}^{l^{[j]}(\mathbf{a}) - 1} \mathbb{1}(a_i = o) \ , \ j = 1, \cdots, k+1 \
\ \Bigg\}
\end{split} 
\end{align} 

We have
\begin{align} 
\begin{split}
\mathcal{C}_{\tilde{n}, \mathbf{z}}^{(4)} = 
		\bigcup_{\mathbf{\check a} \in \check{\mathcal{C}}_{\mathbf{z}}^{(2)}} 
		\bigcup_{\vec{\mathbf{z}} \in \vec{\mathcal{Z}}_{\mathbf{z}, \mathbf{\check a}}^{(1)}} 
		\bigcup_{\mathbf{n}_o \in \mathbb{N}_0^{k+1}} 
		\mathcal{C}_{\tilde{n}, \mathbf{\check a}, \vec{\mathbf{z}}, \mathbf{n}_o}^{(5)}
\end{split} 
\end{align} 
Meanwhile, for any $\mathbf{a} \in \mathcal{C}_{\tilde{n}, \mathbf{\check a}, \vec{\mathbf{z}}, \mathbf{n}_o}^{(5)}$, 
\begin{align} 
\begin{split}
& \mathbb{P} \Big[\mathbf{C} = \mathbf{a} \Big| N = \tilde{n} + \sum_{j=1}^{k+1} n_o^{[j]} \Big] (\beta) = \\
	&\quad \prod_{j=1}^{k+1} \Bigg( \Big( P_{o: \mathcal{A}_1 \cup \{o\} \setminus \{\check{a}^{[1]}, \cdots, \check{a}^{[j-1]}\} } (\beta) \Big)^{n_o^{[j]}} 
		\prod_{a \in \mathcal{A}_1} 
		\Big( P_{a: \mathcal{A}_1 \cup \{o\} \setminus \{\check{a}^{[1]}, \cdots, \check{a}^{[j-1]}\} } (\beta) \Big)^{n_a^{[j]} + \mathbb{1}(a = \check{a}^{[j]})} \Bigg) \\
\end{split} 
\end{align} 
which is a constant of $\beta$. Note that 
\begin{align} 
\begin{split}
|\mathcal{C}_{\tilde{n}, \mathbf{\check a}, \vec{\mathbf{z}}, \mathbf{n}_o}^{(5)}| 
= \prod_{j=1}^{k+1} \frac{(n_o^{[j]} + \tilde{n}^{[j]})!}{(n_o^{[j]})! \prod_{a \in \mathcal{A}_1} (n_a^{[j]})!}
\end{split} 
\end{align}
Thus, 
\begin{align} 
\begin{split}
L_5[\tilde{n}, \mathbf{z}] (\lambda, \beta) 
&= \sum_{\mathbf{\check a} \in \check{\mathcal{C}}_{\mathbf{z}}^{(2)}} 
		\sum_{\vec{\mathbf{z}} \in \vec{\mathcal{Z}}_{\mathbf{z}, \mathbf{\check a}}^{(1)}} 
		\sum_{\mathbf{n}_o \in \mathbb{N}_0^{k+1}} 
		\mathbb{P}\Big[N = \tilde{n} + \sum_{j=1}^{k+1} n_o^{[j]} \Big] (\lambda) \cdot 
		\mathbb{P}\Big[\mathbf{C} \in \mathcal{C}_{\tilde{n}, \mathbf{\check a}, \vec{\mathbf{z}}, \mathbf{n}_o}^{(5)} \Big| 
			N = \tilde{n} + \sum_{j=1}^{k+1} n_o^{[j]} \Big] (\beta) \\
&= \sum_{\mathbf{\check a} \in \check{\mathcal{C}}_{\mathbf{z}}^{(2)}} 
		\sum_{\vec{\mathbf{z}} \in \vec{\mathcal{Z}}_{\mathbf{z}, \mathbf{\check a}}^{(1)}} 
		\sum_{\mathbf{n}_o \in \mathbb{N}_0^{k+1}} \Bigg(
		\frac{(T \lambda)^{\tilde{n} + \sum_{j=1}^{k+1} n_o^{[j]}} \cdot e^{-T \lambda}}{(\tilde{n} + \sum_{j=1}^{k+1} n_o^{[j]})!} \cdot 
		\Big( \prod_{j=1}^{k+1} \frac{(n_o^{[j]} + \tilde{n}^{[j]})!}{(n_o^{[j]})! \prod_{a \in \mathcal{A}_1} (n_a^{[j]})!} \Big) \cdot 
\\
&\quad\quad \prod_{j=1}^{k+1} \Big( (P_{o: \mathcal{A}_1 \cup \{o\} \setminus \{\check{a}^{[1]}, \cdots, \check{a}^{[j-1]}\} } (\beta))^{n_o^{[j]}} 
		\prod_{a \in \mathcal{A}_1} 
		(P_{a: \mathcal{A}_1 \cup \{o\} \setminus \{\check{a}^{[1]}, \cdots, \check{a}^{[j-1]}\} } (\beta) )^{n_a^{[j]} + \mathbb{1}(a = \check{a}^{[j]})} \Big) \Bigg) \\
\end{split} 
\end{align} \\

\subsection{Likelihood of Sales under Attraction Demand Models}

In $L_5[\tilde{n}, \mathbf{z} | N \leq m] (\lambda, \beta)$, the size of the summations grow exponentially fast in $k$, and an exact evaluation can easily become impossible. However, the likelihood function can be simplified if we can choose choice models with nice mathematical structures. In this section, we discuss the attraction demand models, and show how it can help reducing the size of the sums, and make the evaluation of likelihood functions more practical. (Note that the widely used multinomial logit model is one special case of the attraction demand models.) The discussion focuses mainly on $L_5[\tilde{n}, \mathbf{z}] (\lambda, \beta)$, but the same idea applies to  $L_4[\tilde{n}, \mathbf{\tilde a}] (\lambda, \beta)$ as well. \\

Consider a generic attraction demand model 
\begin{align} 
\begin{split}
P_{a:\mathcal{A}_1 \cup \{o\}} (\beta) = \frac{f_a(\beta_a)}{1 + \sum_{a' \in \mathcal{A}_1} f_{a'}(\beta_{a'})} \quad ; \quad 
P_{o:\mathcal{A}_1 \cup \{o\}} (\beta) = \frac{1}{1 + \sum_{a' \in \mathcal{A}_1} f_{a'}(\beta_{a'})} \\
\end{split} 
\end{align} 
With this model, we have 
\begin{align} 
\begin{split}
&\quad \mathbb{P}[\mathbf{C} = \mathbf{a} | N = n] (\beta) 
		= \prod_{i=1}^n P_{a_i:\mathcal{A}_i(a_1, \cdots, a_{i-1}) \cup \{o\}} (\beta) 
		= \frac{\prod_{a \in \mathcal{A}_1} (f_a(\beta_a))^{n_a} }
				{\prod_{j=1}^{k+1} \Big(1 + \sum_{a' \in \mathcal{A}_1 \setminus \{\check{a}^{[1]}, \cdots, \check{a}^{[j-1]}\}} f_{a'}(\beta_{a'})\Big)^{n^{[j]} + \mathbb{1}(j<k+1)} }
\end{split} 
\end{align}
where $n_a, \forall a \in \mathcal{A}_1$, $k$ and $n^{[j]}, j \in \{1,\cdots, k+1\}$ are the values of $N_a, \forall a \in \mathcal{A}_1$, $K$ and $N^{[j]}, j \in \{1,\cdots, K+1\}$ given $\mathbf{C} = \mathbf{a}$. As we can see, now the likelihood of observing $\mathbf{C} = \mathbf{a}$ only depends on the value of these three sets of random variables. \\

Recall that for any $\tilde{N} = \tilde{n}$ and any $\mathbf{Z} = \mathbf{z} = (n_a, \forall a \in \mathcal{A}_1) \in \prod_{a \in \mathcal{A}_1} \{0, \cdots, s_a\}$ such that $\tilde{n} = \sum_{a \in \mathcal{A}_1} n_a$, 
\begin{align} 
\begin{split}
L_5[\tilde{n}, \mathbf{z}] (\lambda, \beta) 
&= \mathbb{P}[E_5(\tilde{n}, \mathbf{z})] (\lambda, \beta) 
= \sum_{\mathbf{a} \in \mathcal{C}_{\tilde{n}, \mathbf{z}}^{(4)}} 
		\mathbb{P}[N = |\mathbf{a}|] (\lambda) \cdot 
		\mathbb{P}[\mathbf{C} = \mathbf{a} | N = |\mathbf{a}|] (\beta) 
\end{split} 
\end{align} 
Meanwhile, 
\begin{align} 
\begin{split}
\check{\mathcal{C}}_{\mathbf{z}}^{(2)} \defi \Big\{ \check{\mathbf{a}} = (\check{a}^{[1]}, \cdots, \check{a}^{[k]}) \in \check{\mathcal{C}}_{k}^{(1)} \ : \ 
		n_{\check{a}^{[j]}} = s_{\check{a}^{[j]}} \ , \ \forall j = 1, \cdots, k \ \Big\}
\end{split} 
\end{align} 
For any $\mathbf{\check a} \in \check{\mathcal{C}}_{\mathbf{z}}^{(2)}$, let
\begin{align} 
\begin{split}
\mathcal{N}_{\tilde{n}, \mathbf{z}, \mathbf{\check a}}^{(1)} \defi \Bigg\{ \mathbf{n} = (n^{[1]}, \cdots, n^{[k+1]}) \in \mathbb{N}_0^{k+1} \ : \
		& \sum_{j'=1}^j (n^{[j']}+1) \geq \sum_{j'=1}^j n_{\check{a}^{[j']}} \ , \ j = 1, \cdots, k  \ ; \
		\sum_{j=1}^{k+1} (n^{[j]} + 1) \geq \tilde{n} + 1 \ \Bigg\}
\end{split} 
\end{align} 
denotes the feasible set of value that $\mathbf{N}$ can take given $\mathbf{Z} = \mathbf{z}$ and $\check{\mathbf{C}} = \check{\mathbf{a}}$, where $k$ is the value of $K$ given $\mathbf{Z} =  \mathbf{z}$. \\

Now, for any $\mathbf{\check a} \in \check{\mathcal{C}}_{\mathbf{z}}^{(2)}$ and any $\mathbf{n} \in \mathcal{N}_{\tilde{n}, \mathbf{z}, \mathbf{\check a}}^{(1)}$, let
\begin{align} 
\begin{split}
\mathcal{C}_{\mathbf{z}, \mathbf{\check a}, \mathbf{n}}^{(7)} \defi \Bigg\{ 
		\mathbf{a} \in \mathcal{C}^{(1)}_{\sum_{j=1}^{k} (n^{[j]}+1) - 1} \ : \ 
		n^{[j]} \defi l^{[j]}(\mathbf{a}) - l^{[j-1]}(\mathbf{a}) - 1 \ \Bigg\}
\end{split} 
\end{align} 
We have
\begin{align} 
\begin{split}
\mathcal{C}_{\tilde{n}, \mathbf{z}}^{(4)} = 
		\bigcup_{\mathbf{\check a} \in \check{\mathcal{C}}_{\mathbf{z}}^{(2)}} 
		\bigcup_{\mathbf{n} \in \mathcal{N}_{\tilde{n}, \mathbf{z}, \mathbf{\check a}}^{(1)}} 
		\mathcal{C}_{\mathbf{z}, \mathbf{\check a}, \mathbf{n}}^{(7)} 
\end{split} 
\end{align} 
With the attraction demand model, for any $\mathbf{a} \in \mathcal{C}_{\mathbf{z}, \mathbf{\check a}, \mathbf{n}}^{(7)}$, the choice likelihood
\begin{align} 
\begin{split}
P \Big[\mathbf{C} = \mathbf{a} | N = \sum_{j=1}^{k} (n^{[j]}+1) - 1 \Big] = 
		\frac{\prod_{a \in \mathcal{A}_1} (f_a(\beta_a))^{n_a} }
			{\prod_{j=1}^{k+1} \Big(1 + \sum_{a' \in \mathcal{A}_1 \setminus \{\check{a}^{[1]}, \cdots, \check{a}^{[j-1]}\}} 
				f_{a'}(\beta_{a'})\Big)^{n^{[j]} + \mathbb{1}(j<k+1)} }
\end{split} 
\end{align}
is a constant of $\beta$. Note that 
\begin{align} 
\begin{split}
|\mathcal{C}_{\mathbf{z}, \mathbf{\check a}, \mathbf{n}}^{(7)}| 
		= \Big( \frac{(\sum_{a \in \mathcal{A}_1 \setminus \{\check{a}^{[1]}, \cdots, \check{a}^{[k]}\}} n_a)!}
			{\prod_{a \in \mathcal{A}_1 \setminus \{\check{a}^{[1]}, \cdots, \check{a}^{[k]}\}} (n_a)!} \Big) 
			\prod_{j=1}^{k} {n^{[j]} + \sum_{j'=1}^{j-1} (n^{[j']} + 1 - n_{\check{a}^{[j]}}) \choose n_{\check{a}^{[j]}}-1}
\end{split} 
\end{align}
Thus, 
\begin{align} 
\begin{split}
L_5[\tilde{n}, \mathbf{z}] (\lambda, \beta) 
&= \sum_{\mathbf{\check a} \in \check{\mathcal{C}}_{\mathbf{z}}^{(2)}} 
		\sum_{\mathbf{n} \in \mathcal{N}_{\tilde{n}, \mathbf{z}, \mathbf{\check a}}^{(1)}} 
		\mathbb{P}\Big[N = \sum_{j=1}^{k} (n^{[j]}+1) - 1 \Big] (\lambda) \cdot 
		\mathbb{P}\Big[\mathbf{C} \in \mathcal{C}_{\tilde{n}, \mathbf{\check a}, \vec{\mathbf{z}}, \mathbf{n}_o}^{(5)} \Big| 
			N = \sum_{j=1}^{k} (n^{[j]}+1) - 1 \Big] (\beta) \\
&= \sum_{\mathbf{\check a} \in \check{\mathcal{C}}_{\mathbf{z}}^{(2)}} 
		\sum_{\mathbf{n} \in \mathcal{N}_{\tilde{n}, \mathbf{z}, \mathbf{\check a}}^{(1)}} \Bigg(
		\frac{(T \lambda)^{\sum_{j=1}^{k} (n^{[j]}+1) - 1} \cdot e^{-T \lambda}}{\Big( \sum_{j=1}^{k} (n^{[j]}+1) - 1 \Big)!} \cdot 
		\Big( \frac{(\sum_{a \in \mathcal{A}_1 \setminus \{\check{a}^{[1]}, \cdots, \check{a}^{[k]}\}} n_a)!}
			{\prod_{a \in \mathcal{A}_1 \setminus \{\check{a}^{[1]}, \cdots, \check{a}^{[k]}\}} (n_a)!} \Big) 
			\prod_{j=1}^{k} {n^{[j]} + \sum_{j'=1}^{j-1} (n^{[j']} + 1 - n_{\check{a}^{[j]}}) \choose n_{\check{a}^{[j]}}-1} \cdot 
\\
&\quad\quad\quad \frac{\prod_{a \in \mathcal{A}_1} (f_a(\beta_a))^{n_a} }
			{\prod_{j=1}^{k+1} \Big(1 + \sum_{a' \in \mathcal{A}_1 \setminus \{\check{a}^{[1]}, \cdots, \check{a}^{[j-1]}\}} 
					f_{a'}(\beta_{a'})\Big)^{n^{[j]} + \mathbb{1}(j<k+1)} } \Bigg) \\
\end{split} 
\end{align} \\

\section{Likelihood of Sales without the Null Alternative} 

In this part, we discuss the likelihood of sales under a different setting: when the null alternative is not available ($s_o = 0$). This is the setting which [Conlon and Mortimer, 2013] works with. Based on the discussion, we will provide a counter-example to show that the formulation of likelihood function given by [Conlon and Mortimer, 2013] is incorrect. \\

Consider the scenario when we observe $\tilde{N} = \tilde{n}$ for some $\tilde{n} \in \mathbb{N}_0$, and $\mathbf{Z} = \mathbf{z}$ for some $\mathbf{z} = (n_a \in \mathcal{A}_1)$, where $n_a \in \{0, \cdots, s_a\}$ for all $a \in \mathcal{A}$. Since $s_o = 0$, we have $N = \tilde{N} = \tilde{n}$. Meanwhile, given $\mathbf{Z} = \mathbf{z}$, the values of $K$ is known. Let it be $k$. Let 
\begin{align} 
\begin{split}
E_6(\tilde{n}, \mathbf{n}) \defi \Big\{ \omega \in \Omega \ : \ N_a(\omega) = n_a \ , \ \forall a \in \mathcal{A}_1 \ ; \ N_o(\omega) = 0 \Big\}
\end{split} 
\end{align} 
denotes the event that we observe $\mathbf{Z} = \mathbf{z}$ and $N_o = 0$. We have
\begin{align} 
\begin{split}
E_6(\tilde{n}, \mathbf{n}) = \Big\{ \omega \in \Omega \ : \ \mathbf{C}(\omega) \in \mathcal{C}_{\tilde{n}, \mathbf{n}}^{(5)} \Big\}
\end{split} 
\end{align} 
where
\begin{align} 
\begin{split}
\mathcal{C}_{\tilde{n}, \mathbf{n}}^{(5)} \defi \Big\{ \mathbf{a} \in \mathcal{C}_{\tilde n}^{(1)} :  
	\quad \sum_{i=1}^{\tilde{n}} \mathbb{1}(a_i = a) = n_a \ , \ \forall a \in \mathcal{A}_1 \ \Big\} 
\end{split} 
\end{align} 
is the set of choice sequences $\mathbf{C}$ that are consistent with $\mathbf{n}$, and contain no choice of the null alternative. \\ 

Let $\mathbb{P}'$ be the modified probability measure, on which we assume that the null alternative is not available. The only way $\mathbb{P}'$ being different from $\mathbb{P}$ is that 
\begin{align} 
\begin{split}
&\quad \mathbb{P}'[\mathbf{C} = \mathbf{a} | N = n] (\beta) 
	= \mathbb{P}'[\mathbf{C} = \mathbf{a} | \mathbf{C} \in \mathcal{C}^{(1)}_{n}] (\beta) \\
	&= \mathbb{P}'[C_1 = a_1] (\beta) \cdot \mathbb{P}'[C_2 = a_2 | C_1 = a_1] (\beta) 
			\cdots \mathbb{P}'[C_n = a_n | C_{n-1} = a_{n-1}, \cdots, C_1 = a_1] (\beta) \\
	&= \prod_{i=1}^n P_{a_i:\mathcal{A}_i(a_1, \cdots, a_{i-1})} (\beta) 
\end{split} \label{measure-seq}
\end{align}
Thus, likelihood of observing $\mathbf{Z} = \mathbf{z}$ is
\begin{align} 
\begin{split}
L_6[\tilde{n}, \mathbf{z}] (\lambda, \beta)
&= \mathbb{P}'[E_6(\tilde{n}, \mathbf{n}) | N=\tilde{n}] (\beta) 
= \mathbb{P}'[N=\tilde{n}] (\lambda) \cdot 
		\sum_{\mathbf{a} \in \mathcal{C}_{\tilde{n}, \mathbf{n}}^{(5)}} \mathbb{P}'[\mathbf{C} = \mathbf{a} | N=\tilde{n}] (\beta) 
\end{split} 
\end{align} 

For any $\check{\mathbf{a}} \in \check{\mathcal{C}}_{\mathbf{z}}^{(2)}$ and any $\vec{\mathbf{z}} \in \vec{\mathcal{Z}}_{\mathbf{z}, \mathbf{\check a}}^{(1)}$, let
\begin{align} 
\begin{split}
\mathcal{C}_{\tilde{n}, \mathbf{\check a}, \vec{\mathbf{z}}}^{(6)} \defi \Bigg\{ 
		\mathbf{a} \in\mathcal{C}^{(1)}_{\tilde{n}} \ : \ 
		\check{a}^{[j]} &= a_{l^{[j]}(\mathbf{a})} \ , \ j = 1, \cdots, k \ ;  \\
		n_a^{[j]} &= \sum_{i=l^{[j-1]}(\mathbf{a})+1}^{l^{[j]}(\mathbf{a}) - 1} \mathbb{1}(a_i = a) \ , \ j = 1, \cdots, k+1 \ , \ \forall a \in \mathcal{A}_1 
\ \Bigg\}
\end{split} 
\end{align} 
We have
\begin{align} 
\begin{split}
\mathcal{C}_{\tilde{n}, \mathbf{n}}^{(5)} = 
		\bigcup_{\mathbf{\check a} \in \check{\mathcal{C}}_{\mathbf{z}}^{(2)}} 
		\bigcup_{\vec{\mathbf{z}} \in \vec{\mathcal{Z}}_{\mathbf{z}, \mathbf{\check a}}^{(1)}} 
		\mathcal{C}_{\tilde{n}, \mathbf{\check a}, \vec{\mathbf{z}}}^{(6)}
\end{split} 
\end{align} 
Meanwhile, for any $\mathbf{a} \in \mathcal{C}_{\tilde{n}, \mathbf{\check a}, \vec{\mathbf{z}}}^{(6)}$, 
\begin{align} 
\begin{split}
& \mathbb{P}'[\mathbf{C} = \mathbf{a} | N=\tilde{n}] (\beta) = 
		\prod_{j=1}^{k+1} \prod_{a \in \mathcal{A}_1} 
		\Big( P_{a: \mathcal{A}_1 \setminus \{\check{a}^{[1]}, \cdots, \check{a}^{[j-1]}\} } (\beta) \Big)^{n_a^{[j]} + \mathbb{1}(a = \check{a}^{[j]})}  \\
\end{split} 
\end{align} 
which is a constant of $\beta$. Note that 
\begin{align} 
\begin{split}
|\mathcal{C}_{\tilde{n}, \mathbf{\check a}, \vec{\mathbf{z}}}^{(6)}| 
= \prod_{j=1}^{k+1} \frac{(n^{[j]})!}{\prod_{a \in \mathcal{A}_1} (n_a^{[j]})!}
\end{split} 
\end{align}
Thus, 
\begin{align} 
\begin{split}
& L_6[\tilde{n}, \mathbf{z}] (\lambda, \beta) \\
&=  \mathbb{P}'[N=\tilde{n}] (\lambda) \cdot 
		\sum_{\mathbf{\check a} \in \check{\mathcal{C}}_{\mathbf{z}}^{(2)}} 
		\sum_{\vec{\mathbf{z}} \in \vec{\mathcal{Z}}_{\mathbf{z}, \mathbf{\check a}}^{(1)}} 
		\mathbb{P}'[\mathbf{C} \in \mathcal{C}_{\tilde{n}, \mathbf{\check a}, \vec{\mathbf{z}}}^{(6)} | N=\tilde{n}] (\beta) \\
&= \frac{(T \lambda)^{\tilde{n}} \cdot e^{-T \lambda}}{(\tilde{n})!} \cdot 
		\sum_{\mathbf{\check a} \in \check{\mathcal{C}}_{\mathbf{z}}^{(2)}} 
		\sum_{\vec{\mathbf{z}} \in \vec{\mathcal{Z}}_{\mathbf{z}, \mathbf{\check a}}^{(1)}} \Bigg(
		\Big( \prod_{j=1}^{k+1} \frac{(n^{[j]})!}{\prod_{a \in \mathcal{A}_1} (n_a^{[j]})!} \Big) \cdot 
		\prod_{j=1}^{k+1} \prod_{a \in \mathcal{A}_1} 
		\Big( P_{a: \mathcal{A}_1 \setminus \{\check{a}^{[1]}, \cdots, \check{a}^{[j-1]}\} } (\beta) \Big)^{n_a^{[j]} + \mathbb{1}(a = \check{a}^{[j]})} \Bigg) \\
\end{split} 
\end{align} 

Under attraction demand models, the likelihood function becomes
\begin{align} 
\begin{split}
L_6[\tilde{n}, \mathbf{z}] (\lambda, \beta) 
&= \frac{(T \lambda)^{\tilde{n}} \cdot e^{-T \lambda}}{(\tilde{n})!} \cdot 
		\Big( \frac{(\sum_{a \in \mathcal{A}_1 \setminus \{\check{a}^{[1]}, \cdots, \check{a}^{[k]}\}} n_a)!}
			{\prod_{a \in \mathcal{A}_1 \setminus \{\check{a}^{[1]}, \cdots, \check{a}^{[k]}\}} (n_a)!} \Big) \cdot 
\\
&\quad\quad\quad 
		\sum_{\mathbf{\check a} \in \check{\mathcal{C}}_{\mathbf{z}}^{(2)}} 
		\sum_{\mathbf{n} \in \mathcal{N}_{\tilde{n}, \mathbf{z}, \mathbf{\check a}}^{(1)}} 
		\Bigg( \prod_{j=1}^{k} {n^{[j]} + \sum_{j'=1}^{j-1} (n^{[j']} + 1 - n_{\check{a}^{[j]}}) \choose n_{\check{a}^{[j]}}-1} \cdot 
			\frac{\prod_{a \in \mathcal{A}_1} (f_a(\beta_a))^{n_a} }
			{\prod_{j=1}^{k+1} \Big(1 + \sum_{a' \in \mathcal{A}_1 \setminus \{\check{a}^{[1]}, \cdots, \check{a}^{[j-1]}\}} 
					f_{a'}(\beta_{a'})\Big)^{n^{[j]} + \mathbb{1}(j<k+1)} } \Bigg) \\
\end{split} 
\end{align} 

Recall that the complete data likelihood function $L_1$ and $L_2$ can be separated into two parts, with one part only depends on $\lambda$, and the other only depends on $\beta$. Although the same argument does not apply to $L_3$, $L_4$ and $L_5$, it does apply to $L_6$. As we can see, since the null alternative is not available, $\tilde{n}^{[j]} = n^{[j]}$ for any $j = 1, \cdots, k+1$. Thus, $\tilde{n} = (\sum_{j=1}^{k+1} (n^{[j]}+1)) - 1$ always holds true, which allows us to separate the arrival likelihood $(T \lambda)^{\tilde{n}} \cdot e^{-T \lambda} / (\tilde{n})!$. Let
\begin{align} 
\begin{split}
L_6^{(2)}[\tilde{n}, \mathbf{z}] (\beta) 
&= \sum_{\mathbf{\check a} \in \check{\mathcal{C}}_{\mathbf{z}}^{(2)}} 
		\sum_{\vec{\mathbf{z}} \in \vec{\mathcal{Z}}_{\mathbf{z}, \mathbf{\check a}}^{(1)}} \Bigg(
		\Big( \prod_{j=1}^{k+1} \frac{(n^{[j]})!}{\prod_{a \in \mathcal{A}_1} (n_a^{[j]})!} \Big) \cdot 
		\prod_{j=1}^{k+1} \prod_{a \in \mathcal{A}_1} 
		\Big( P_{a: \mathcal{A}_1 \setminus \{\check{a}^{[1]}, \cdots, \check{a}^{[j-1]}\} } (\beta) \Big)^{n_a^{[j]} + \mathbb{1}(a = \check{a}^{[j]})} \Bigg) 
\end{split} 
\end{align} 
be the second part of $L_6[\tilde{n}, \mathbf{z}] (\lambda, \beta)$. The likelihood function [Conlon and Mortimer, 2013] tries to derive is $L_6^{(2)}[\tilde{n}, \mathbf{z}] (\beta) $, which corresponds to the conditional choice probability $\mathbb{P}'[\mathbf{C} \in \mathcal{C}_{\tilde{n}, \mathbf{n}}^{(5)} | N=\tilde{n}] (\beta)$, and only depends on $\beta$.

\subsection{Some Counter-examples}

Except some discussions in the appendix, [Conlon and Mortimer, 2013] focuses mainly on the case when there is only one stock-out event. Let $a^*$ denotes the only product that is out of stock by time $T$. Let $N^* \defi \sum_{a \in \mathcal{A}_1 \setminus \{a^*\}} N_a^{[1]}$ denotes the total number of customers that purchase products other than $a^*$ before $a^*$ becomes out of stock. \\

[Conlon and Mortimer, 2013] first argues that $N^*$ follows a truncated negative binomial distribution between $0$ and $n - n_{a^*}$. However, this is not true in general. To see this, simply assume that there is a product $a'$, which is "dominated" by $a^*$ (e.g. similar product but with slightly different quality). That is, $a'$ can never be chosen when $a^*$ is available, but has a positive probability to be chosen when $a^*$ is out of stock. When $n_{a'}>0$, the probability that $N^*$ takes value of $n - n_{a^*}$ is zero. (Note that we are discussing the conditional distribution. If we the total number of arrivals are not observed, then $N^*$ follows a negative binomial distribution.) \\
 
The argument does hold true, however, when we have
\begin{align} 
\begin{split}
P_{a : A} (\beta) = (1 - P_{a^* : \{a, a^*\}} (\beta)) \cdot P_{a : A \setminus \{a^*\}} (\beta) \ , \ \forall a \neq a^* \ , \ a \in \mathcal{A}
\end{split} 
\end{align}
In this case, the likelihood function is reduced to
\begin{align} 
\begin{split}
&\quad L_6[\tilde{n}, \mathbf{z}] (\lambda, \beta) \\
&= \Big( \frac{(\tilde{n} - n_{a^*})!}{\prod_{a \in \mathcal{A}_1 \setminus \{a^*\}} (n_a)!} \cdot 
		\prod_{a \in \mathcal{A}_1 \setminus \{a^*\}} (P_{a:\mathcal{A}_1 \setminus \{a^*\}} (\beta))^{n_a} \Big) 
		\sum_{n^*=0}^{\tilde{n} - n_{a^*}}  \Bigg( \Big( \frac{(n_{a^*} - 1 + n^*)!}{(n_{a^*} - 1)! (n^*)!} \Big) \cdot 
		( P_{a^* : \{a, a^*\}} (\beta) )^{n_{a^*}} \cdot (1 - P_{a^* : \mathcal{A}_1} (\beta))^{n^*} \Bigg) \\ 
\end{split} \label{mle6-1}
\end{align}
where $N^*$ follows a truncated negative binomial distribution between $0$ and $n - n_{a^*}$ with success rate $P_{a^* : \mathcal{A}_1}$. \\

[Conlon and Mortimer, 2013] further argues that given $N^*$, each individual sale $N_a^{[1]}$ follows a binomial distribution independently for every $a \in \mathcal{A}$, while the parameter of this distribution depends on $N^*$. In any case, the fact that $N^* = \sum_{a \in \mathcal{A}_1 \setminus \{a^*\}} N_a^{[1]}$ means the distributions of $N_a^{[1]}$ cannot be independent. Thus, the second argument is not correct. \\

To make it more concrete, consider the case when $\mathcal{A} = \{a, a^*\}$, where $a^*$ is out of stock, and $a$ is still available by time $T$. Following the previous setup in this section, we assume that we observe $N_{a^*} = n_{a^*}$ for some $n_{a^*} \in \mathbb{N}$, and $N_a = n_a$ for some $n_a \in \mathbb{N}$. From here, we can make a clear comparison of the two different formulations. \\ 

[Conlon and Mortimer, 2013] does not give the likelihood function explicitly. Instead, it suggests the use of EM algorithm, which maximizes
\begin{align} 
\begin{split}
n_{a^*} \log(P_{a^* : \{a, a^*\}}) + E[N_a^{[1]}] \log(P_{a : \{a, a^*\}}) + E[n_a - N_a^{[1]}] \log(P_{a : \{a\}})
\end{split} 
\end{align}
Note that $\log(P_{a : \{a\}}) = \log(1) = 0$. Thus, to prove our argument, we can simply compare the value of $E[N_a^{[1]}]$ obtained from the two models. Recall that the distribution of $N^*$ is
\begin{align} 
\begin{split}
p(n^*) \defi \mathbb{P}(N^* = n^*) = \frac{\Big( \frac{(n_{a^*} - 1 + n^*)!}{(n_{a^*} - 1)! (n^*)!} \Big) \cdot 
		( P_{a^* : \{a, a^*\}} (\beta) )^{n_{a^*}} \cdot (1 - P_{a^* : \{a, a^*\}} (\beta))^{n^*}}
		{\sum_{n'=0}^{\tilde{n} - n_{a^*}}  \Big( \Big( \frac{(n_{a^*} - 1 + n')!}{(n_{a^*} - 1)! (n')!} \Big) \cdot 
		( P_{a^* : \{a, a^*\}} (\beta) )^{n_{a^*}} \cdot (1 - P_{a^* : \{a, a^*\}} (\beta))^{n'} \Big)} \\
\end{split} 
\end{align}
[Conlon and Mortimer, 2013] argues that given $N^* = n^*$, $N_a^{[1]} \sim \text{Binomial}(n_a, \rho(n^*))$, where
\begin{align} 
\begin{split}
\rho(n^*) &= \frac{n^* P_{a : \{a, a^*\}}}{n^* P_{a : \{a, a^*\}} + (n_a - n^*) P_{a : \{a\}}} 
	= \frac{n^* P_{a : \{a, a^*\}}}{n_a - n^* (1 - P_{a : \{a, a^*\}}) }
\end{split} 
\end{align}
Thus, it is suggesting
\begin{align} 
\begin{split}
E'[N_a^{[1]}] = \frac{\sum_{n^*=0}^{n_a} \Big( 
			n_a \cdot \frac{n^* P_{a : \{a, a^*\}}}{n_a - n^* (1 - P_{a : \{a, a^*\}}) } \cdot
			\Big( \frac{(n_{a^*} - 1 + n^*)!}{(n_{a^*} - 1)! (n^*)!} \Big) \cdot 
			( P_{a^* : \{a, a^*\}} (\beta) )^{n_{a^*}} \cdot (1 - P_{a^* : \{a, a^*\}} (\beta))^{n^*} \Big)}
		{\sum_{n'=0}^{\tilde{n} - n_{a^*}}  \Big( \Big( \frac{(n_{a^*} - 1 + n')!}{(n_{a^*} - 1)! (n')!} \Big) \cdot 
		( P_{a^* : \{a, a^*\}} (\beta) )^{n_{a^*}} \cdot (1 - P_{a^* : \{a, a^*\}} (\beta))^{n'} \Big)} 
\end{split} \label{em-cm2}
\end{align} 
However, in this scenario, it is very clear that $N_a^{[1]} = n^*$. Thus by (\ref{mle6-1}), 
\begin{align} 
\begin{split}
E[N_a^{[1]}] = E[n^*] = \frac{\sum_{n^*=0}^{n_a} \Big(n^* \cdot \Big( \frac{(n_{a^*} - 1 + n^*)!}{(n_{a^*} - 1)! (n^*)!} \Big) \cdot 
		( P_{a^* : \{a, a^*\}} (\beta) )^{n_{a^*}} \cdot (1 - P_{a^* : \{a, a^*\}} (\beta))^{n^*} \Big)}
		{\sum_{n'=0}^{\tilde{n} - n_{a^*}}  \Big( \Big( \frac{(n_{a^*} - 1 + n')!}{(n_{a^*} - 1)! (n')!} \Big) \cdot 
		( P_{a^* : \{a, a^*\}} (\beta) )^{n_{a^*}} \cdot (1 - P_{a^* : \{a, a^*\}} (\beta))^{n'} \Big)}
\end{split} \label{em-correct}
\end{align}
We can see that $E[N_a^{[1]}] \neq E'[N_a^{[1]}]$. As an example, let $n_{a^*} = n_a = 2$, and let $P_{a : \{a, a^*\}} = 1/2$. Then (\ref{em-correct}) suggests that $E[N_a^{[1]}] = 10/11$, while (\ref{em-cm2}) suggests that $E'[N_a^{[1]}]=26/33$. This shows the formulation in [Conlon and Mortimer, 2013] is not correct. \\

\section{Estimation}

$L_1[n, \mathbf{a}, \mathbf{t}] (\lambda, \beta)$ and $L_2[n, \mathbf{a}] (\lambda, \beta)$ are easy to work with, since the logarithms of both functions are concave in $\lambda$ and $\beta$. $L_3[\tilde{n}, \mathbf{\tilde a}, \mathbf{\tilde t}] (\lambda, \beta) (\lambda, \beta)$ is not concave in general, but the function is still easy to evaluate. Meanwhile, for any fixed value of $\lambda$, the logarithm of $L_3[\tilde{n}, \mathbf{\tilde a}, \mathbf{\tilde t}] (\lambda, \beta) (\lambda, \beta)$ is concave in $\beta$. Since $\lambda$ is one-dimensional, we can perform a line search over $\lambda$ to find its optimal value. When it comes to estimating the parameters from $L_4[\tilde{n}, \mathbf{\tilde a}] (\lambda, \beta)$ or $L_5[\tilde{n}, \mathbf{z}] (\lambda, \beta)$, however, things can be much harder. \\

In this section, we discuss two ways that can make the evaluation of the likelihood functions more practical. The discussion focuses on $L_5[\tilde{n}, \mathbf{z}] (\lambda, \beta)$, while similar ideas can be applied to $L_4[\tilde{n}, \mathbf{\tilde a}] (\lambda, \beta)$ as well. \\

\subsection{Dealing with Infinite Sums}

In the previous sections, we discussed the use of attraction demand models, and showed that it can reduce the size of the sums. However, the sums are still infinite. In practice, we can find an integer $m$ that is large enough, such that the likelihood $\mathbb{P}[N > m] (\lambda)$ is small enough to be ignored. In this case, we can approximate $L_5[\tilde{n}, \mathbf{z}] (\lambda, \beta)$ with the conditional likelihood function $L_5[\tilde{n}, \mathbf{z} | N \leq m] (\lambda, \beta)$. Since the Poisson likelihood decays quickly, usually $m$ doesn't need to be very large. \\

The calculation of $L_5[\tilde{n}, \mathbf{z} | N \leq m] (\lambda, \beta)$ only involves finite sum:
\begin{align} 
\begin{split}
L_5[\tilde{n}, \mathbf{z} | N \leq m] (\lambda, \beta) 
&= \sum_{\mathbf{\check a} \in \check{\mathcal{C}}_{\mathbf{z}}^{(2)}} 
		\sum_{\vec{\mathbf{z}} \in \vec{\mathcal{Z}}_{\mathbf{z}, \mathbf{\check a}}^{(1)}} \ \
		\sum_{\mathbf{n}_o \in \mathbb{N}_0^{k+1}, \ n_o^{[1]} + \cdots + n_o^{[k+1]} \leq \tilde{n} - m} \Bigg(
		\frac{(T \lambda)^{\tilde{n} + \sum_{j=1}^{k+1} n_o^{[j]}} \cdot e^{-T \lambda}}{(\tilde{n} + \sum_{j=1}^{k+1} n_o^{[j]})!} \cdot 
		\Big( \prod_{j=1}^{k+1} \frac{(n_o^{[j]} + \tilde{n}^{[j]})!}{(n_o^{[j]})! \prod_{a \in \mathcal{A}_1} (n_a^{[j]})!} \Big) \cdot 
\\
&\quad\quad \prod_{j=1}^{k+1} \Big( (P_{o: \mathcal{A}_1 \cup \{o\} \setminus \{\check{a}^{[1]}, \cdots, \check{a}^{[j-1]}\} } (\beta))^{n_o^{[j]}} 
		\prod_{a \in \mathcal{A}_1} 
		(P_{a: \mathcal{A}_1 \cup \{o\} \setminus \{\check{a}^{[1]}, \cdots, \check{a}^{[j-1]}\} } (\beta) )^{n_a^{[j]} + \mathbb{1}(a = \check{a}^{[j]})} \Big) \Bigg) 
\end{split} 
\end{align} 

To derive $L_5[\tilde{n}, \mathbf{z} | N \leq m] (\lambda, \beta)$ under attraction demand models, let 
\begin{align} 
\begin{split}
\mathcal{N}_{n, \mathbf{z}, \mathbf{\check a}}^{(2)} \defi \Bigg\{ \mathbf{n} = (n^{[1]}, \cdots, n^{[k+1]}) \in \mathbb{N}^{k+1} \ : \
		& \sum_{j'=1}^j n^{[j']}+1 \geq \sum_{j'=1}^j n_{\check{a}^{[j']}} \ , \ j = 1, \cdots, k  \ ; \
		\sum_{j=1}^{k+1} (n^{[j]} + 1) = n + 1 \Bigg\}
\end{split} 
\end{align} 
We have
\begin{align} 
\begin{split}
\mathcal{C}_{\tilde{n}, \mathbf{z}}^{(4)} 
= \bigcup_{\mathbf{\check a} \in \check{\mathcal{C}}_{\mathbf{z}}^{(2)}} 
		\bigcup_{\mathbf{n} \in \mathcal{N}_{\tilde{n}, \mathbf{z}, \mathbf{\check a}}^{(1)}} 
		\mathcal{C}_{\mathbf{z}, \mathbf{\check a}, \mathbf{n}}^{(7)} 
= \bigcup_{n=\tilde{n}}^{\infty} \
		\bigcup_{\mathbf{\check a} \in \check{\mathcal{C}}_{\mathbf{z}}^{(2)}} 
		\bigcup_{\mathbf{n} \in \mathcal{N}_{n, \mathbf{z}, \mathbf{\check a}}^{(2)}} 
		\mathcal{C}_{\mathbf{z}, \mathbf{\check a}, \mathbf{n}}^{(7)} 
\end{split} 
\end{align} 
Therefore, 
\begin{align} 
\begin{split}
L_5[\tilde{n}, \mathbf{z} | N \leq m] (\lambda, \beta) 
&= \sum_{n=\tilde{n}}^{m} \ \sum_{\mathbf{\check a} \in \check{\mathcal{C}}_{\mathbf{z}}^{(2)}} 
		\sum_{\mathbf{n} \in \mathcal{N}_{n, \mathbf{z}, \mathbf{\check a}}^{(2)}} \Bigg(
		\frac{(T \lambda)^{\sum_{j=1}^{k} (n^{[j]}+1) - 1} \cdot e^{-T \lambda}}{\Big( \sum_{j=1}^{k} (n^{[j]}+1) - 1 \Big)!} \cdot 
		\Big( \frac{(\sum_{a \in \mathcal{A}_1 \setminus \{\check{a}^{[1]}, \cdots, \check{a}^{[k]}\}} n_a)!}
			{\prod_{a \in \mathcal{A}_1 \setminus \{\check{a}^{[1]}, \cdots, \check{a}^{[k]}\}} (n_a)!} \Big) \cdot 
\\
&\quad\quad\quad \prod_{j=1}^{k} {n^{[j]} + \sum_{j'=1}^{j-1} (n^{[j']} + 1 - n_{\check{a}^{[j]}}) \choose n_{\check{a}^{[j]}}-1} \cdot 
			\frac{\prod_{a \in \mathcal{A}_1} (f_a(\beta_a))^{n_a} }
			{\prod_{j=1}^{k+1} \Big(1 + \sum_{a' \in \mathcal{A}_1 \setminus \{\check{a}^{[1]}, \cdots, \check{a}^{[j-1]}\}} 
					f_{a'}(\beta_{a'})\Big)^{n^{[j]} + \mathbb{1}(j<k+1)} } \Bigg) \\
\end{split} \label{cond-mle5-attrac}
\end{align} \\

\subsection{Vectors of Stock-out Indices}

Consider ($\ref{cond-mle5-attrac}$). Although $L_5[\tilde{n}, \mathbf{z} | N \leq m] (\lambda, \beta)$ is greatly simplified under attraction demand models, the size of the set $\bigcup_{\mathbf{\check a} \in \check{\mathcal{C}}_{\mathbf{z}}^{(2)}} \mathcal{N}_{n, \mathbf{z}, \mathbf{\check a}}^{(2)}$ can still be too large. In this case, we can use sampling methods to evaluate $L_5[\tilde{n}, \mathbf{z} | N \leq m] (\lambda, \beta)$. To do that, we need to know the size of the set $\bigcup_{\mathbf{\check a} \in \check{\mathcal{C}}_{\mathbf{z}}^{(2)}} \mathcal{N}_{n, \mathbf{z}, \mathbf{\check a}}^{(2)}$, and to be able to sample from that set  uniformly. In this part, we introduce such a sampling method that we found easy to implement. \\

Consider any $h \in \mathbb{N}$ and any $\check{a}_1', \cdots, \check{a}_{h}' \in \mathcal{A}_1$, such that for any $i \neq j \ , \ i, j \in \{1,\cdots,h\}$, we have $\check{a}_i' \neq \check{a}_j'$. Suppose that we have an arrival sequence with length $n'$, where $n' \geq s_{\check{a}_1'} + \cdots + s_{\check{a}_{h}'}$. During this sequence, the products $\check{a}_1', \cdots, \check{a}_{h}'$ are stocked out. Let $r_{\check{a}_1'}, \cdots, r_{\check{a}_{h}'}$ denote the indices in the arrival sequence where $\check{a}_1', \cdots, \check{a}_{h}'$ become out of stock. We call the vector of those stock-out Indices $(r_{\check{a}_1'}, \cdots, r_{\check{a}_h'})$ a stock-out vector. \\

Clearly, we need $1 \leq r_{\check{a}_j'} \leq n'$ for any $j \in \{1,\cdots,h\}$. Meanwhile, to make $(r_{\check{a}_1'}, \cdots, r_{\check{a}_h'})$ a feasible stock-out vector, the following two rules need to hold as well: 
\begin{align} 
\begin{split}
r_{\check{a}_i'} &\neq r_{\check{a}_j'} \ , \ i \neq j \ , \ i, j \in \{1,\cdots,h\} \\
r_{\check{a}_j'} &\geq \sum_{i=1}^h \mathbb{1}(r_{\check{a}_i'} < r_{\check{a}_j'}) \cdot s_{\check{a}_i'}  
\end{split} 
\end{align}
Let $f(\check{a}_1', \cdots \check{a}_h', n')$ be the total number of distinct, feasible stock-out vectors $(r_{\check{a}_1'}, \cdots, r_{\check{a}_h'})$ give $n'$and $\check{a}_1', \cdots \check{a}_h'$. We have
\begin{align} 
\begin{split}
f(\check{a}_1', \cdots \check{a}_h', n') &= 
		\frac{(n'+1)!}{(n'+1-h)!} - \frac{(n')!}{(n'+1-h)!} \cdot (s_{\check{a}_1'} + \cdots + s_{\check{a}_h'})
\end{split} 
\end{align}
We can prove this equation by induction. First, when $h=1$, the number of distinct vectors $(r_{\check{a}_1'})$ is $n' + 1 - s_{\check{a}_1'}$, which is consistent with the equation above. Second, suppose that for some $h \in \{1, \cdots, |\mathcal{A}|-1\}$, the equation above holds true. Consider the case when there are $h+1$ stock-out events, and let $\check{a}_{h+1}'$ be that new "added" product that is out of stock. Note that any of those $h+1$ products $\check{a}_1', \cdots, \check{a}_{h+1}'$ can be the last one that is out of stock. Meanwhile, the other $h$ stock-outs events must appear before the last one. Thus, conditioning on the arrival index of the last stock-out event and the corresponding product, we have 
\begin{align} 
\begin{split}
f(\check{a}_1', \cdots, \check{a}_{j+1}', n') 
		&= \sum_{j=1}^h \ \sum_{r=s_{\check{a}_1'} + \cdots + s_{\check{a}_{h+1}'}}^{n'} 
			\Bigg( \frac{(r)!}{(r-h)!} - \frac{(r-1)!}{(r-h)!} \cdot ((s_{\check{a}_1'} + \cdots + s_{\check{a}_{h+1}'}) - s_{\check{a}_j'}) \Bigg) \\
		&= \Bigg( \sum_{r=s_{\check{a}_1'} + \cdots + s_{\check{a}_{h+1}'}}^{n'} 
			\frac{(r)!(h+1)}{(r-h)!} \Bigg)
		- \Bigg( \sum_{r=s_{\check{a}_1'} + \cdots + s_{\check{a}_{h+1}'}}^{n'} 
			\frac{(r-1)!(h)}{(r-h)!} \cdot (s_{\check{a}_1'} + \cdots + s_{\check{a}_{h+1}'}) \Bigg) \\
		&= (h+1)! \sum_{r=s_{\check{a}_1'} + \cdots + s_{\check{a}_{h+1}'}}^{n'} {r \choose h}
			-  (s_{\check{a}_1'} + \cdots + s_{\check{a}_{h+1}'}) \cdot (h)! 
				\sum_{r=s_{\check{a}_1'} + \cdots + s_{\check{a}_{h+1}'}}^{n'} {r-1 \choose h-1} 
\end{split} 
\end{align}
Note that
\begin{align} 
\begin{split}
\sum_{r=s_{\check{a}_1'} + \cdots + s_{\check{a}_{h+1}'}}^{n'} {r \choose h} 
	&= \sum_{r=h+1}^{n'} {r \choose h} \ - \ \sum_{r=h+1}^{s_{\check{a}_1'} + \cdots + s_{\check{a}_{h+1}'} - 1} {r \choose h} 
	= {n'+1 \choose h+1} - {s_{\check{a}_1'} + \cdots + s_{\check{a}_{h+1}'} \choose h+1} \\
\sum_{r=s_{\check{a}_1'} + \cdots + s_{\check{a}_{h+1}'}}^{n'} {r-1 \choose h-1} 
	&= \sum_{r=h}^{n'} {r-1 \choose h-1} \ - \ \sum_{r=h}^{s_{\check{a}_1'} + \cdots + s_{\check{a}_{h+1}'} - 1} {r-1 \choose h-1} 
	= {n' \choose h} - {s_{\check{a}_1'} + \cdots + s_{\check{a}_{h+1}'}-1 \choose h} \\
\end{split} 
\end{align}
Therefore, 
\begin{align} 
\begin{split}
f(\check{a}_1', \cdots, \check{a}_{h+1}', n') 
		&= \frac{(n'+1)!}{(n'-h)!} - \frac{(s_{\check{a}_1'} + \cdots + s_{\check{a}_{h+1}'})!}{(n'-h)!} \\
		&- (s_{\check{a}_1'} + \cdots + s_{\check{a}_{h+1}'}) \cdot \frac{(n')!}{(n'-h)!} 
			+ (s_{\check{a}_1'} + \cdots + s_{\check{a}_{h+1}'}) \cdot \frac{(s_{\check{a}_1'} + \cdots + s_{\check{a}_{h+1}'} - 1)!}{(n'-h)!} \\
		&= \frac{(n'+1)!}{(n'+1-(h+1))!} - \frac{(n')!}{(n'+1-(h+1))!} \cdot (s_{\check{a}_1'} + \cdots + s_{\check{a}_{h+1}'})
\end{split} 
\end{align}
which proves that the equation holds for the case of $h+1$. Enough to conclude. \\

Now, come back and consider the scenario when we observe $\tilde{N} = \tilde{n}$ and $\mathbf{Z} = \mathbf{z} = (n_a : a \in \mathcal{A}_1)$. Given $\mathbf{Z} = \mathbf{z}$, let $k$ be the value of $K$, and let $\{\check{a}_1, \cdots, \check{a}_k\}$ be the set of products that are out of stock by time $T$ (in an arbitrary order). Let $r_{\check{a}_1}, \cdots, r_{\check{a}_k}$ denote the indices in the arrival sequence where $\check{a}_1, \cdots, \check{a}_k$ become out of stock, and let $(r_{\check{a}_1'}, \cdots, r_{\check{a}_k})$ denotes the stock-out vector. As we can see, for each $n \geq \tilde{n}$, the set of all feasible stock-out vectors $(r_{\check{a}_1}, \cdots, r_{\check{a}_k})$ has a one-to-one mapping to the set $\bigcup_{\mathbf{\check a} \in \check{\mathcal{C}}_{\mathbf{z}}^{(2)}} \mathcal{N}_{n, \mathbf{z}, \mathbf{\check a}}^{(2)}$: 
\begin{align} 
\begin{split}
\mathbf{\check a} &= (\check{a}^{[1]}, \cdots, \check{a}^{[k]}) \ : \ 
		\{\check{a}^{[1]}, \cdots, \check{a}^{[k]}\} = \{\check{a}_1, \cdots, \check{a}_k\} \ ; \  
		r_{\check{a}^{[1]}} < \cdots < r_{\check{a}^{[k]}} 
\\
n^{[1]} &= r_{\check{a}^{[1]}} - 1 \\
n^{[j]} &= r_{\check{a}^{[j]}} - r_{\check{a}^{[j-1]}} - 1 \ , \ j = 2, \cdots, k \\
n^{[k+1]} &= n - r_{\check{a}^{[k]}}
\end{split} 
\end{align}

Recall that $r_{\check{a}_j} \in \{1, \cdots, n\}$ for $j = 1, \cdots, k$. That is, for each $n$, there are at most $n^k$ distinct vectors of stock-out indices. Therefore, to sample elements from $\bigcup_{\mathbf{\check a} \in \check{\mathcal{C}}_{\mathbf{z}}^{(2)}} \mathcal{N}_{n, \mathbf{z}, \mathbf{\check a}}^{(2)}$ uniformly, we can go through the following steps: (1) Sample integers from $0$ to $n^k - 1$ uniformly (without replacement) using a linear congruential generator. (2) Divide each integer by $n$ for $k$ times, collect the $k$ remainders, and add all of them by $1$ as stock-out indices. (3) Keep only the vectors that are feasible, and transform them into combinations of $\mathbf{\check a}$ and $\mathbf{n}$. \\

Now, we have
\begin{align} 
\begin{split}
\Big| \bigcup_{\mathbf{\check a} \in \check{\mathcal{C}}_{\mathbf{z}}^{(2)}} 
		\mathcal{N}_{n, \mathbf{z}, \mathbf{\check a}}^{(2)} \Big| 
&= f(\check{a}_1, \cdots, \check{a}_k, n) 
	= \frac{(n)!}{(n-k+1)!} \cdot (n + 1 - (s_{\check{a}_1} + \cdots + s_{\check{a}_k})) \\
\end{split} 
\end{align}
Since the total number of distinct integers we can sample from is $n^k$, the acceptance ratio of samples is
\begin{align} 
\begin{split}
\frac{\Big| \bigcup_{\mathbf{\check a} \in \check{\mathcal{C}}_{\mathbf{z}}^{(2)}} 
		\mathcal{N}_{n, \mathbf{z}, \mathbf{\check a}}^{(2)} \Big|}{n^k}
&= \frac{(n)}{n} \cdot \frac{(n-1)}{n} \cdot \cdots \cdot \frac{(n-k+2)}{n} \cdot \frac{n + 1 - (s_{\check{a}_1} + \cdots + s_{\check{a}_k})}{n} \\
\end{split} 
\end{align}
In practice, $n$ is usually much larger than $k$ and $s_{\check{a}_1} + \cdots + s_{\check{a}_k}$ (otherwise it will hurt customers' satisfaction). Therefore, the acceptance ratio is usually large, which makes this method applicable. \\



\section{Numerical Tests and Results}

Consider five products: 0, 1, 2, 3 and 4. Suppose that customer choices follow an attraction demand model
\begin{align} 
\begin{split}
P_{a: \mathcal{A}}(\beta) &= \frac{\beta_{a}}{\sum_{a' \in \mathcal{A}} \beta_{a'} } 
\end{split}
\end{align} 
where $\beta_0 = 0.25$, $\beta_1 = 0.05$, $\beta_1 = 0.1$, $\beta_1 = 0.2$ and $\beta_4 = 0.4$. Product $0$ is always available, and can never go out-of-stock (stock level is infinite). Product $1$, $2$, $3$ and $4$ is offered with probability $0.6$ during each time period (called "visit"). If offered, the initial stock level is $3$. When the sales of a product reaches its stock level, it becomes unavailable to be chosen. Customers arrive according to a Poisson process, with rate $6$ per visit. We simulate $10,000$ such visits. The initial assortment and the accumulative sales of each product is recorded. \\

\begin{figure}[h!]
	\centering
		\includegraphics[width=1\columnwidth]{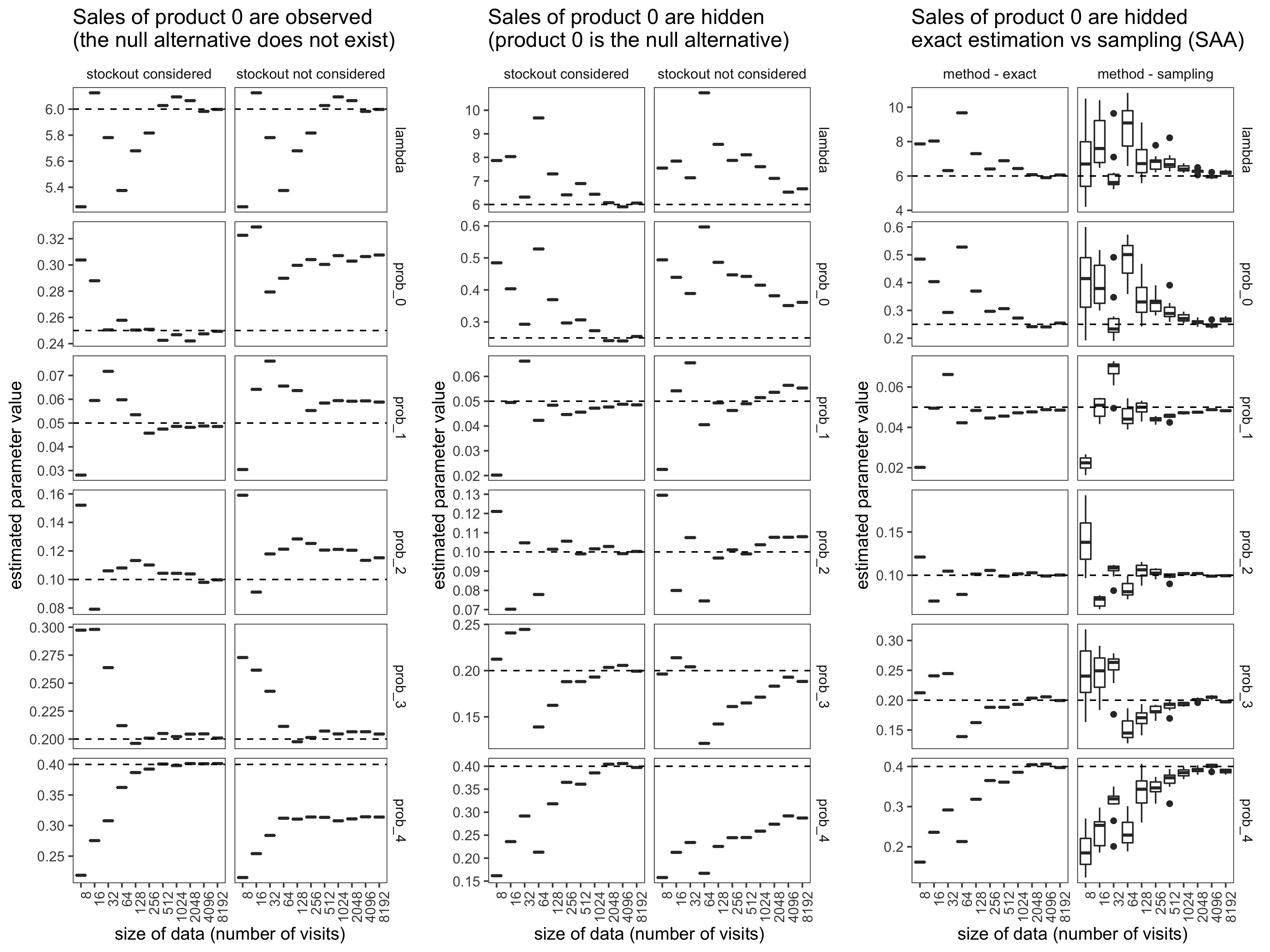} 
		\caption{Results of Numerical Tests}
	\label{plot}
\end{figure}

Based on this simulated data, we performed three tests. The results are presented in Figure \ref{plot}. Each column of plots in the figure corresponds to estimations from a specific likelihood function, while each row of plots corresponds to a specific parameter. The y-axis of each plot shows the values of estimations, the x-axis shows the sizes of data (number of visits) used in those estimations, while the dashed line represents the true values of the parameters. Note that in the model, we can add all $\beta$ parameters by any constant without changing the choice probabilities. Thus, we present the estimated choice probabilities of products instead of the $\beta$ parameters, in order to make the comparison more consistent. \\

The results of the first test is presented in the first and second columns of plots in Figure \ref{plot}. In this test, we estimate two likelihood functions with the sales summary data observed. Column 1 are estimations obtained with $L_6$, which is the correct likelihood function in this case. Column 2 are estimations from the classical likelihood function which ignores stock-outs. \\ 

The results of the second test is presented in the third and forth columns of plots in Figure \ref{plot}. Similarly, we estimate two likelihood functions. However, we hide the sales of product 0 in the entire dataset. In this case, product 0 represents the null alternative whose sales are hidden. Column 3 are estimations obtained with $L_5$, which is the correct likelihood function in this case. Column 4 are estimations from the classical likelihood function which ignores stock-outs and only modeled the null alternative. \\ 

The results of the third test is presented in the fifth and sixth columns of plots in Figure \ref{plot}. In this test, we also hide the sales of product 0 in the entire dataset. Column 5 is identical to Column 3, showing estimations obtained with $L_5$. Column 6 shows the estimations obtained from $L_5$ using the sampling (SAA) method, instead of exact evaluations. The "boxes" give the quantiles from 10 trails. In each trail, we sample one stock-out vector for each possible number of arrivals in each visit, and estimate the parameter values from the approximated likelihood function. \\

\begin{figure}[h!]
	\centering
		\includegraphics[width=1\columnwidth]{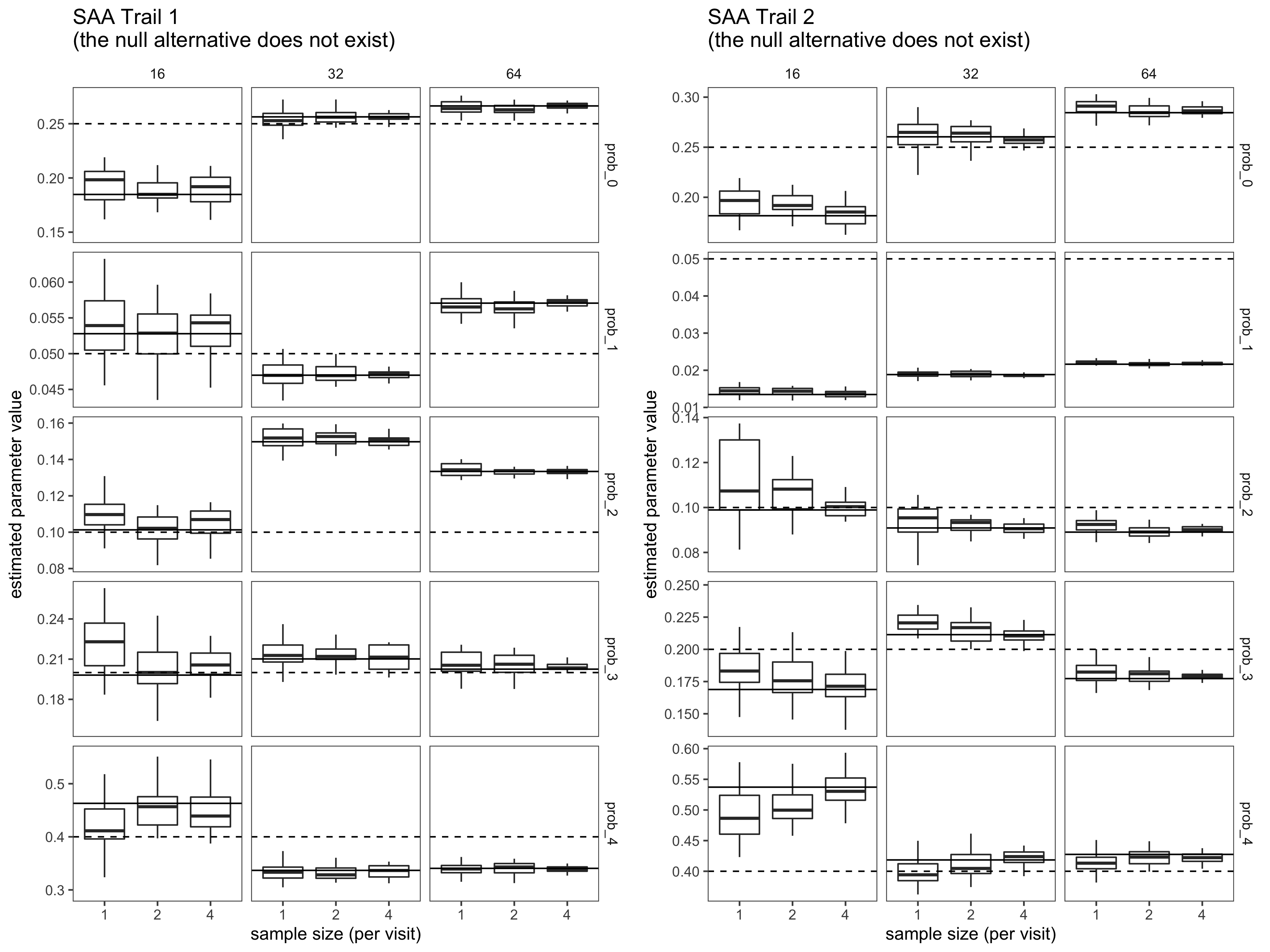} 
		\caption{Results of Numerical Tests}
	\label{plot2}
\end{figure}

We also did an additional test to illustrate how the sampling method works. In this test, we estimate $L_5$ with different amount of data (number of visits), and different sample size (number of stock-out vectors sampled without replacement) per visit. Under each scenario, we perform 20 estimations, and present the results as boxplots (showing the minimum, maximum and the quantiles). We repeat the text two times with different simulated data. All results from this test are presented in Figure \ref{plot2}. Each column of plots in the figure corresponds to estimations obtained from a specific size of data, while each row of plots corresponds to a specific parameter. The y-axis of each plot shows the values of estimations, the x-axis shows the sample size (number of stock-out vectors per visit) used in those estimations, while the solid line represents the estimated values obtained from exact calculations. The true value of the parameters are still represented by dashed lines. Clearly, estimations obtained from the sampling method gets better with more data, and more samples as well. \\

As a short conclusion, these results show that ignoring stock-outs can lead to serious bias in estimation. In the third test, we can see that we a reasonable amount of data, we can get good estimations even with a very small sample size. It shows that our model is practical enough to be used. \\

\section{Conclusion}

As the numerical results showed, if unobserved stock-out events are not modeled appropriately, they can introduce bias to the estimation of discrete choice models. With the likelihood function being correctly formulated, the bias can be corrected. On the other hand, although likelihood functions with stock-out events can be much harder to evaluate, it is still possible to get good estimations quickly using stochastic optimization methods. \\

There are still questions remain unanswered. Among them, one important research opportunity lies in the identifiability of the likelihood functions. For complete data and transaction data with transaction times, researchers have already developed the necessary and sufficient conditions for the parameters to be identifiable from the likelihood function. For example, the arrival rate is identifiable from the likelihood of transaction data with transaction times, if and only if there is a sufficient variety in choice sets. In other scenarios, however, identification conditions remain unknown: Since stock-out events introduce more variety of choice sets, the number of potential customers may be identifiable, even if the initial assortment is a constant over each time period. More work is needed to address such problems. \\

\newpage

\bibliography{reference}{}
\bibliographystyle{apalike}

\newpage
\section{Appendix: }

\end{document}